\documentclass[12pt,twoside]{article}
\usepackage{caption}
\usepackage{keyval}
\usepackage{graphicx} 
\usepackage{url}
\usepackage{booktabs}
\usepackage{color,xcolor}
\usepackage[utf8]{inputenc}
\usepackage[T1]{fontenc}
\usepackage[a4paper, left=2.5cm, right=2.5cm, top=3cm, bottom=3cm]{geometry}
\usepackage{mathptmx} 
\usepackage{pgfgantt}
\usepackage{tikz}
\usepackage{enumitem}
\setlist{nosep}
\usepackage[sectionbib]{natbib}

\usepackage{hyperref}
\usepackage{lscape}
\usepackage{color} 
\usepackage{sectsty}
\usepackage{chngcntr}

\usepackage{blindtext}
\usepackage{fancyhdr}
\usepackage{amsmath,amssymb}

\newtheorem{notation}{Notation}[section]

\newcommand{\cL}{\mathcal{L^*}}
\newcommand{\Pool}{\mathcal{L^{\textup{Pool}}}}
\newcommand{\PTN}{G}
\newcommand{\NN}{\mathbb{N}}

\newcommand{\trips}{trips}
\newcommand{\kante}{a}
\newcommand{\ODset}{\mathcal{O}\hspace{-0.7mm}\mathcal{D}\hspace{-0.5mm}}
\newcommand{\OD}{\textup{OD}}
\newcommand{\odpair}{(v^{\textrm{org}},v^{\textrm{dest}})}
\newcommand{\stations}{V}
\newcommand{\org}{v^{\textrm{org}}}
\newcommand{\dest}{v^{\textrm{dest}}}
\newcommand{\Pcgn}{P^{\textrm{CGN}}}
\newcommand{\Pean}{P^{\textrm{EAN}}}
\newcommand{\Pptn}{P^{\textrm{PTN}}}
\newcommand{\cgn}{\textrm{CGN}}
\newcommand{\cgnnodeset}{N}

\newcommand{\cgnarcset}{E}

\newcommand{\paxnumber}{q}

\newcommand{\cost}{c}
\newcommand{\costvariable}{\textrm{cost}^v}
\newcommand{\cU}{{\cal U}}
\title{
Planning and optimizing transit lines
}
\author{
Marie Schmidt University Würzburg, Würzburg, Germany
\\
marie.schmidt@uni-wuerzburg.de\\[1ex]
Anita Schöbel, University Kaiserslautern-Landau and\\
Fraunhofer ITWM, Kaiserslautern, Germany
\\
anita.schoebel@mathematik.rptu.de
\\
}

\begin{document}
\maketitle

\paragraph{Abstract}

For all line-based transit systems like bus, metro and tram, the routes
of the lines and the frequencies at which they are operated are determining
for the operational performance of the system. However,
as transit line planning happens early in the planning process, it is not straightforward
to predict the effects of line planning decisions on relevant performance indicators.
This challenge has in more than 40 years of research on transit line planning
let to many different models.
In this chapter, we concentrate on models for transit line planning
including transit line planning under uncertainty.
We pay particular attention to the interplay 
of passenger routes, frequency and capacity, and specify three different
levels of aggregation at which these can be modeled. 

Transit line planning has been studied in different communities under different names.
The problem can be decomposed into the components \emph{line generation}, \emph{line selection}, and \emph{frequency setting}. We include publications that
regard one of these individual steps as well as publications that combine
two or all of them.
We do not restrict to models build with a certain solution approach in mind,
but do have a focus on models expressed in the language of mathematical programming.



\section{Introduction}

\begin{figure}
\begin{center}
\includegraphics[width=0.5 \textwidth]{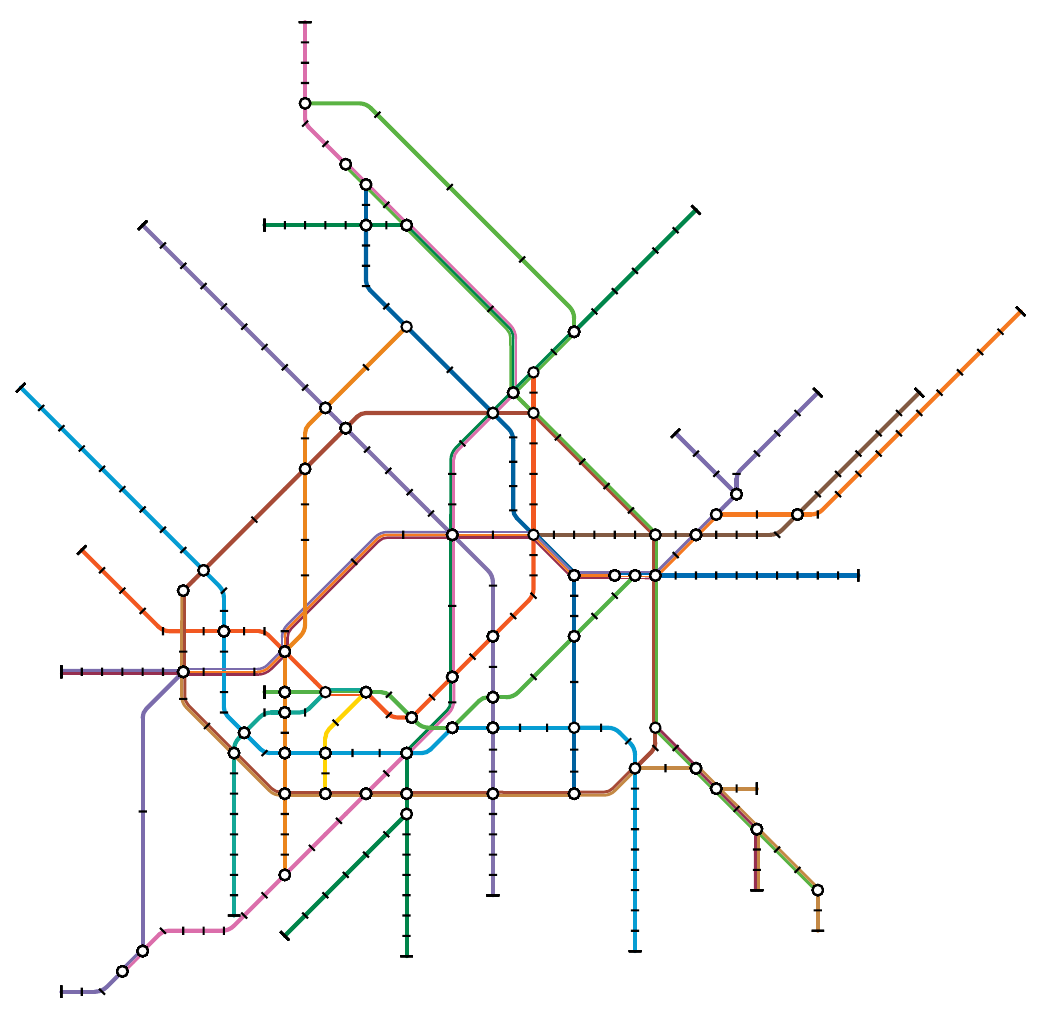}
\end{center}
\caption{A possible visualization of the line plan of the city of Berlin 
  \label{fig-lineplan2}}
\end{figure}

\paragraph{What is transit line planning?}
Colloquially speaking, transit line planning is to decide which public transport
lines are to be established when the stops or stations and possible direct connections between them are already known, and demand between these can be estimated.
 More precisely, the decisions
concern the number of lines to be established, their routes and the frequencies
with which they should be operated. The result is a line plan together with
the frequencies of the lines. 
The lines themselves can be depicted on a map,
as, e.g., the well-known map of the London underground. As an example, a
map of the line plan of the city of Berlin (including \emph{U-Bahn}
and \emph{S-Bahn} lines operational in 2009) is shown in Figure~\ref{fig-lineplan2}.

\paragraph{Relevance of transit line planning.}
Designing a network of transit lines and deciding upon
their frequencies is a core process in transit planning.
The decision on which locations are connected by direct lines can have a big impact on mobility patterns in the considered region.
In schedule-based transportation, transit line planning sets
the ground for the subsequent steps such as timetabling, vehicle- and crew
scheduling, delay management and others. 
As such, transit line planning 
is considered the basis for the public transport supply and 
the decisions taken at this planning stage as fundamental for 
most performance indicators (such as efficiency, emissions, costs).
In particular in a metropolitan area, a transit network coordinates
different modes such as buses, light rail, trams, or metros.

Nevertheless, line planning is a relevant problem also for railway applications.
We hence also refer to some papers which are formulated in the railway context.

\paragraph{Different names - clarification.}
Transit line planning has mostly been researched in \emph{transport engineering} and in
  \emph{operations research}. This explains a hodgepodge of subproblems and variations with
  different names and notations.
  In order to give a systematic view, we split transit line planning into
  three subproblems (see Figure~\ref{fig-clarification}) showing the
  overlap of TNDP and line planning. The three subproblems of
  \emph{transit line planning} are:

\begin{figure}
\begin{center}
\includegraphics[width=0.8 \textwidth]{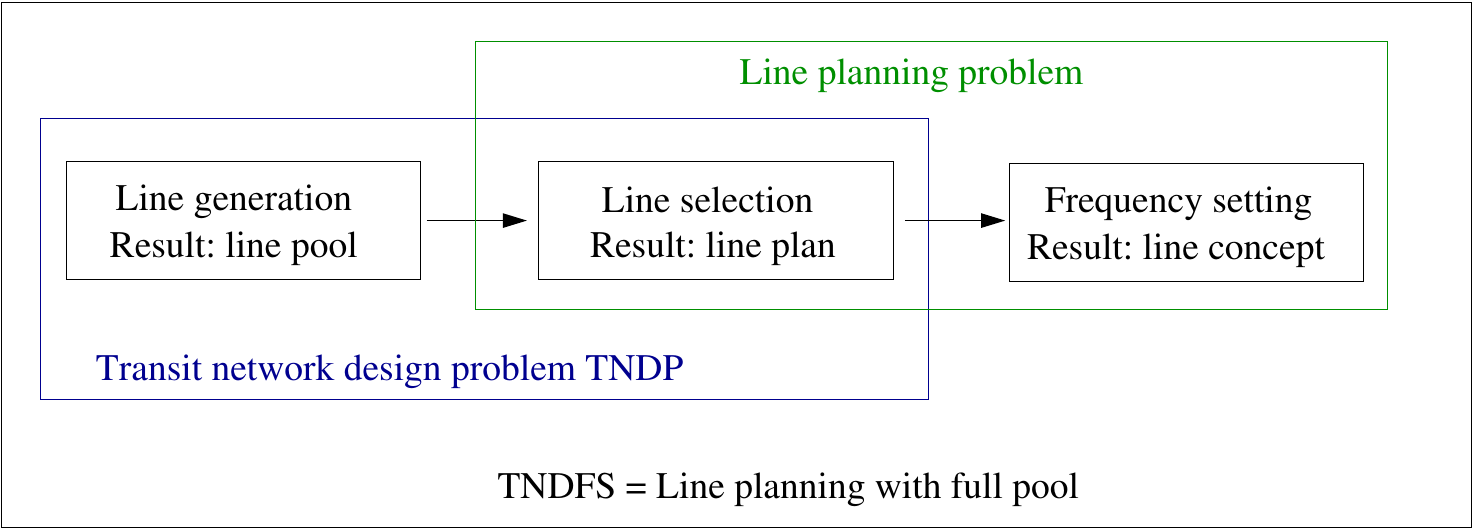}
\end{center}
\caption{Subproblems of transit line planning and their combinations}
\label{fig-clarification}
\end{figure}  

\begin{enumerate}
\item \emph{Line generation}: a set of candidate lines is constructed.
  The result is the \emph{line pool}. 
\item \emph{Line selection}: From a given line pool, the lines to be established
  are selected. The result is the \emph{line plan}.
\item \emph{Frequency setting}: Frequencies are assigned to the selected lines.
  The lines together with the frequencies are called \emph{line concept}.
\end{enumerate}  
\smallskip

The \emph{transit (route) network design problem (TNDP)} (also called (TRNDP))
mostly researched in transport engineering, covers the subproblems
line generation and line selection.
A line pool as intermediate result is not needed in the TNDP.
Frequency setting is left as separate step after the lines have
been decided. The problem that integrates line generation, line selection, and frequency
setting is known as \emph{transit network design and frequency setting (TNDFS)} in the
transport engineering community, see, e.g., \cite{pternea2015sustainable,jha2019multi}.

The term \emph{line planning}, mostly used in operations research,
integrates the subproblems line selection and frequency setting.
As input, most researchers assume a line pool as given.
There are also publications in which the lines are constructed from scratch within
the line planning problem, e.g., \cite{BGP07,heinrich2022algorithms}. Such models
integrating line generation, line selection and frequency setting
are in operations research called \emph{line planning with the full pool},
since \emph{all} possible lines can be selected in this setting.

\emph{Line planning with the full pool} and \emph{transit network design and frequency
setting (TNDFS)} are hence synonyms for the integration of all three subproblems.
In this paper we use the notation \emph{transit line planning} whenever 
line generation, line selection, frequency setting or a combination of them
is involved.

\paragraph{Development of transit line planning.}
The first paper on transit line planning we are aware of is \cite{Patz25},
published a century ago.
Although dealing with the same type of problem, there has been little overlap
between the literature in the transport engineering and  operations research communities
in the early years.
Nevertheless, both communities mention 
only single papers from the 60s, while research on the topic started in the 70s and evolved
to numerous papers in the 90s and in the first decade of this century.
Since then, research on transit line planning
has continued and the overlap between the communities has continuously increased.

Several surveys exist, most of them from more than a decade ago: 
The survey \cite{Sch10b} is mainly on line planning (line selection and frequency setting)
while the surveys \cite{CedWil86,GuiHao08,KepKar09,FaraMia13,CancMautUru15} are mainly
about the transit network design problem (TNDP).
The recent survey by \cite{duran2022survey} concentrates on developments in TNDP und TNDFS
since 2009.
In this chapter we concentrate on modeling aspects of transit line planning
including transit line planning under uncertainty. Moreover, we specify
three aggregation levels and classify approaches with respect to them.
\medskip

In contrast to vehicle and crew scheduling, algorithms for transit line planning are
still not commonly used in commercial software systems, as mentioned, e.g., in
  \cite{duran2022survey,Migl2024-heureka}. This is about to be changed.
  In the open-source research library {\sf LinTim}, see \cite{lintimhp}, many algorithms and data
sets for transit line planning are available and the method of 
\cite{gatt_et_al:OASIcs.ATMOS.2022.5} is currently further developed for the use in
the professional Heur\`es software of
\cite{heures}.

\begin{figure}
\begin{center}
\includegraphics[width=0.9 \textwidth]{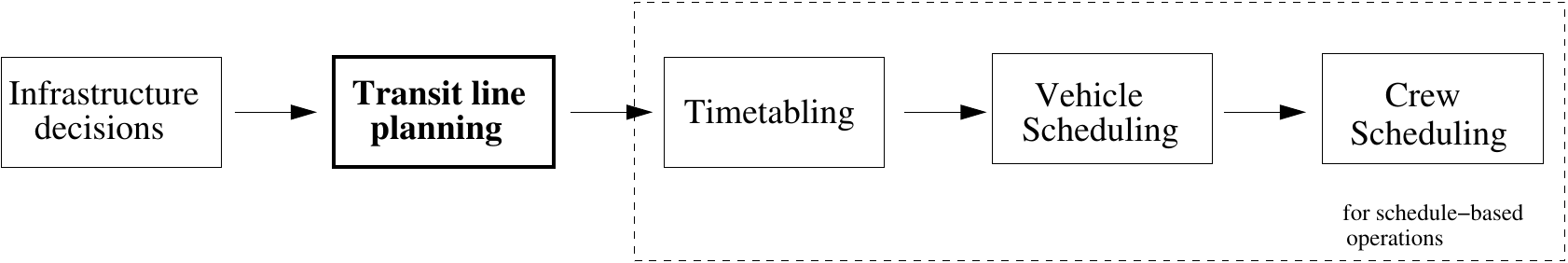} 
\end{center}
\caption{The planning process in public transportation. \label{fig-stages}}
\end{figure}

\paragraph{Embedding of transit line planning in the planning horizon.}\label{embedding}
The planning process in transit is normally divided into several stages where the output of a stage serves as input to the subsequent stage,
see Figure~\ref{fig-stages} which can be found in similar form in many papers
from the 1980s to now, see, e.g., \cite{CedWil86,LiTang23}.

The first stage within the planning horizon,
\emph{infrastructure location decisions}, decides about stops and their direct connections.
 It is followed by \emph{transit line planning}, which is the focus of this chapter.
 In \emph{schedule-based} operations, the next planning stages are \emph{timetabling}
specifying the exact departure and arrival times for the transit lines, and
\emph{vehicle scheduling} assigning vehicles to \trips. In a \emph{line-pure schedule},
each vehicle is assigned to one line which it serves back and forth the whole day. \emph{Crew scheduling} is the next stage.
In \emph{headway-based} 
transportation (also referred to as \emph{schedule-free} transportation),
operations are not meant to adhere to fixed timetables,
but only the headway of a line is known, i.e., the time between two subsequent departures.
Passengers hence go to a stop and wait there for their bus coming.

Other planning steps include crew rostering \cite{CapTothVigoFisch98} or fare
setting \cite{OttoBoysen17,HaSc01}.

\paragraph{Solution approaches for transit line planning.}
Transit line planning invites a variety of solution approaches, spanning from
\emph{construction heuristics}, 
\emph{metaheuristics},
and \emph{mixed-integer programming} approaches to \emph{game theoretic} approaches.  Since this chapter is mainly about \emph{models},
we recommend \cite{KepKar09} to learn more about transit line planning \emph{heuristics},
including early contributions.
The review papers \cite{iliopoulou2019metaheuristics,duran2022survey} give a good overview on
transit line planning with
metaheuristics. For a game theoretic approach we refer to \cite{SchieweSchieweSchmidt19}.

There is also a stream of research aiming to derive optimal frequencies for transit lines analytically by continuum approximation, starting with the research by \cite{mohring1972optimization} on a single line
and later been extended to spatially diversified demand, compare, e.g.,
\cite{jara2003single,fielbaum2016optimal}. This class of models with
parametric descriptions of network structure and demand is described in
Section~\ref{section-parametric}.
\bigskip

\paragraph{Remainder of this chapter.}
Section~\ref{sec-modelingbasic} describes the basic notation, introduces a first
line planning model and discusses modeling of the fundamental elements:
line pool and frequencies. The section also covers performance indicators for
transit line planning.
Passengers play an important role in line planning. Their
routes relate the three key concepts \emph{frequency}, \emph{capacity},
and \emph{demand}. This relation, and the many modeling approaches for routing 
passengers within transit line planning are discussed and classified in
Section~\ref{sec-pax-and-cap}.
Section~\ref{sec-uncertainty} describes how uncertainty can be incorporated by
robust and stochastic transit line planning models. We treat uncertainty of
demand, of driving times and link failures.
Extensions and related problems are discussed in Section~\ref{sec-extensions}. This includes
the skip-stop problem, seasonal demand, the integration with other planning stages
and parametric transit line planning as well as a sketch of other interesting related problems.
We conclude the chapter in Section~\ref{sec-conclusion}.

\section{Basic modeling considerations for transit line planning}\label{sec-modelingbasic}

In the following we first introduce the basic notation needed and specify the variables
and the general goals of transit line planning in Section~\ref{sec-basic}. With
this notation it is already possible to present the basic feasibility constraints and
a first model as a building block for possible extensions
in Section~\ref{sec-costmodel}. Modeling the line pool
is treated in Section~\ref{sec-linepool} and different ways of modeling
the frequencies are shown in Section~\ref{sec-freq}.

\subsection{Basic notation, variables and general goals}\label{sec-basic}

Many different versions of the transit line planning problem have been formulated in the
literature. So there is not a single 'transit line planning model', but there are many
models for each of the subproblems line generation, line selection and frequency
setting or for their combinations.
They differ in the level of accuracy in which they model important aspects like demand or cost, but also in the objective function and in their constraints.

To formally define a line, we first need a model for the underlying infrastructure.

\begin{notation}\label{nota-PTN}
The \emph{public transportation network (PTN)}
$\PTN=(\stations,A)$ consists of nodes $\stations$ which represent stops or stations
and of links (also called \emph{arcs}) $A$ which represent direct (i.e., non-stop)
connections between stops. Link labels in the PTN represent distance which can be
measured as physical distance in meters or as travel time in minutes. 
\end{notation}

A \emph{line} corresponds to a path in the public transport network $\PTN$.
In the basic transit line planning models, it is assumed that the line stops at every
node of the path.
The line generation problem constructs a \emph{line pool} $\Pool$ of candidate lines.
The goal of the line selection problem is to select lines 
to be operated from the line pool $\Pool$. We refer to the subset of selected
lines $\cL\subset \Pool$ as
\emph{line plan}. Sometimes, binary \emph{indicator variables}
$y_l \in \{0,1\}$, $l \in \Pool$ are used to model which lines are selected.

The frequency $f_l$ of a line $l$ denotes how often the line is operated per
time period $T$ (often: per hour) and is usually required to be a natural number.
The frequency of a line has a crucial impact on cost (km driven, number of
vehicles needed, number of crew members needed) and service level (expected waiting
time, capacity per hour). 
For a given line plan $\cL$, the \emph{frequency setting problem} determines
frequencies to operate the already selected lines $l\in \cL$.
Let $(f_l)_{l \in \cL}$ be the vector
which contains the frequencies of all lines $l \in \cL$.
A line plan $(\cL, (f_l)_{l \in \cL})$ together with its vector of frequencies is called a
\emph{line concept}. 

The \emph{line planning problem} (see Figure~\ref{fig-clarification})
asks for a line concept, i.e., for both the lines and their frequencies.
These two subproblems can be modeled integratedly by using the frequency
vector $(f_l)_{l \in \Pool}$ for \emph{all} lines $l \in \Pool$.
$f_l=0$ means that line $l \in \Pool$ is not selected, i.e.,
\[ \cL:=\{l \in \Pool: f_l >0\}. \]
The variables $y_l$ may not be needed in this integrated setting, but can be obtained
by using $y_l=1$ if and only if $f_l >0$.
\bigskip

For evaluating a line concept, performance indicators
can be split into two groups: 
\begin{itemize}
\item \emph{demand-centered indicators} can be coverage, 
  travel time, directness, or the number of connections per hour.
  Their general tendency is that  'more is better', i.e., increasing lines, frequency, or capacity normally has a positive effect on these indicators.

In order to evaluate demand-related indicators, we assume that
    for every pair $(s,t)$ with $s,t \in V$ it is known how many passengers
    $\OD_{s,t}$ wish to travel from station $s$ to station $t$.
    The corresponding matrix $\OD=\left( \OD_{st} \right)_{(s,t) \in V \times V}$ is called
    \emph{OD-matrix}. $\ODset=\{(s,t)\in V \times V : \OD_{s,t} >0\}$ denotes the set of
    origin-destination (OD) pairs.
  
  \item \emph{supply-centered indicators} include emissions and different types of
    cost, e.g., related to distance driven, or vehicles or crew needed. Here,
    the heuristic principle  'less is better' holds, as every additional line or frequency
    increases cost and emissions.

For computing supply-centered indicators we assume that costs or emissions can
    be estimated. This includes costs/emissions per kilometer traveled, dependent on the
    vehicle type used, and costs per hour traveled reflecting, e.g., cost for the driver
    or conductor. 
\end{itemize}

All transit line planning models need to find a balance between these two groups of indicators
and integrate them at least in a basic way. 
\smallskip

It is important to note that the evaluation of a line plan can only be a (rough)
estimate of the actual realization of cost and quality:
As long as only the line plan and the frequencies of the lines are known, the cost and
the travel times cannot be computed exactly, but only estimated since travel time
depends on the timetable and cost depends on the vehicle- and crew schedules.
This gives room for different performance indicators described in more detail in Section~\ref{sec-performance-indicators}.

\subsection{A basic line planning model}\label{sec-costmodel}

We start with one of the simplest integer linear programming models from the literature, compare, e.g.,  \cite{Sch10b}. It solves the line selection and frequency setting problems integratedly (with a given line pool) and hence belongs to the class of line planning models.

\paragraph{Basic feasibility constraint.}
The following basic feasibility constraints are part of most transit line planning models. Let $\Pool$ be the given line pool. As decision variables we use the frequencies
$f_l$ for $l\in \Pool$. We require them to be non-negative and integer.
A frequency of $f_l=0$ means that line $l$ is not operated. The basic feasibility constraints are
\begin{eqnarray}
    L_{\kante}\ \ \le \sum_{l \in \Pool: l\ni a} f_l \ \ \le \ \ U_{\kante}  &&  \forall {\kante}\in A \label{constraint-costmodel}\\
  f_l \in \NN_0 \ \ && \forall l\in \Pool \label{var-costmodel}.
\end{eqnarray}

The constraints \eqref{constraint-costmodel}  bound the cumulative frequency of
each link $a\in A$ from above and below. 
Upper bounds support the supply-centered indicators. Small upper bounds $U_a$
lead typically to less expensive line concepts. The upper bounds are often also due to the capacity of the underlying infrastructure. In particular in rail-based systems where safety headways need to be respected, there is a strict upper bound on the number of \trips\ that can cross the arc $\kante \in A$ in one period. But there may also be reasons to constrain the frequencies in road-based transport, e.g., noise protection, or avoiding road damage by over usage.  

The lower bounds are due to the demand-centered 
indicators. They impose a
minimum cumulative frequency per link,  motivated by service level 
or by capacity considerations. 

Finding frequencies that fulfill constraints \eqref{constraint-costmodel} and
\eqref{var-costmodel} already is strongly NP-hard, see \cite{Bussieck98}
In the reduction, maximum frequencies of the constructed instance are set to $1$, thus, the proof implies strong NP-hardness of the line selection problem as well.

\paragraph{Basic cost model.}
The basic feasibility constraints can be extended to the following   line planning model which is often used as reference or to show new developments.

Let $\OD_{s,t}$ specify the number of travelers per period between $s$ and
$t$, $s,t \in V$ in the PTN. In the cost model, the
lower bounds $L_\kante$ for each link $\kante$ are determined
as follows: Each OD-pair $\OD_{s,t}$
is routed along a shortest path $P_{s,t}$ in the PTN. For every link $\kante$,
  \begin{equation}
    \label{trafficload}
    d_\kante:=\sum_{(s,t): \kante \in P_{s,t}} \OD_{s,t} 
  \end{equation}
  is the number of travelers (the \emph{traffic load}) on link $\kante$.
  Let $C$ be the (constant) capacity of the vehicles. Then,
  $L_\kante:= \lceil \frac{d_\kante}{C} \rceil$ 
is the minimum number of vehicles needed per period along
  link $\kante$ to ensure that all passengers can travel on a shortest path $P_{s,t}$.

Let $c_l$ be the cost of line $l \in \Pool$. It is often assumed that it is
  composed by a constant fixed cost and costs related to the length of the line
  and the time needed for a complete round-trip.
  Then the \emph{basic cost model} 
  \begin{eqnarray}
    \label{costmodel}
    \min \textup{cost}((\Pool,f)) & := & \sum_{l \in \Pool} c_l \cdot f_l\\
    \mbox{s.t. } \ \ L_{\kante} & \le  & \sum_{l \in \Pool: l\ni a} f_l \ \ \le \ \ U_{\kante} \ \ \forall {\kante}\in A \nonumber \\
    f_l & \in  & \NN_0 \ \ \forall l\in \mathcal{L} \nonumber
   \end{eqnarray} 
  finds a line concept with minimal cost in which all passengers can travel on
  their shortest paths in the PTN. The model originally stems from \cite{CDZ98} and
  has been used in many other papers, e.g.,
  \cite{TorTorBoPf08,Sch10b,FHSS17,csahin2020multi,multiperiod}.

\subsection{Line pool generation}\label{sec-linepool}

  The basic cost model \eqref{costmodel} is formulated with a given line pool as input.
  However, the choice of a suitable line pool
  to use as input for the line selection problem is non-trivial. In fact, the choice of a good line pool has a significant impact on the quality of the line concepts,
  as shown in \cite{Heureka17,Migl2024-heureka}

  There are, to the extent of our knowledge, only few papers that study line pool generation in isolation.
  \cite{GHS16} show that the line generation problem is (strongly) NP-hard if 
  a set with a bounded number of lines satisfying constraints \eqref{constraint-costmodel}
  is searched.

  There is a number of publications that study line pool generation together with the subsequent transit line planning steps, i.e., the (TNDP), in a sequential manner.
 In a first step, promising lines are created as paths in a PTN, in a second step, lines are selected. The result of the first step then is a line pool. This has 
been done already in the early papers by
\cite{silman1974planning,CedWil86}, but is still used in more recent contributions, e.g., \cite{suman2019improvement,VerEngPhiVansteen21,gioialine}. 

In numerous contributions that use metaheuristics for solving the TNDP,
one or several line plans are proposed, evaluated,
and subsequently improved in each iteration, see \cite{iliopoulou2019metaheuristics,duran2022survey} for an overview.

The idea of generating promising lines during the optimization process
is also gaining popularity in the operations research community, for example,
by column generation in \cite{BGP07,bertsimas2021data}.
 \cite{heinrich2022algorithms,heinrich_et_al:OASIcs.ATMOS.2023.4}  prove complexity results and state exact  algorithms for the TNDP on specific graph structures.
 
 Most problem formulations 
 for line generation work under the assumption that a public transport system, and in  particular the  line pool, needs to be constructed \emph{from scratch}.  A notable exception is \cite{suman2019improvement}, where preexisting lines can be replaced by new lines with same start and end terminals under the constraint that line length may not deviate by more than a given factor from the shortest path from start to end terminal. Having generated all lines that fulfill these properties, they aim at choosing a subset of them that maximizes the number of additional direct travelers.
This is modeled as a linear integer program. 
\cite{friedrich2021kombination} argue that existing, manually constructed line plans often contain knowledge of transport planners that should not be discarded lightly, and propose to use this expertise to construct a line pool, from which lines are selected by optimization methods. 
\cite{Migl2024-heureka} propose to generate line pools by combining automated steps
and steps that include manual choices.

\subsection{Modeling frequencies}\label{sec-freq}

The basic feasibility constraints \eqref{constraint-costmodel}, \eqref{var-costmodel}
model frequencies by an integer variable $f_l$ that represents
the number of trips (or line runs) per period.

  \citet{FHSS-CASPT18} argue that to improve memorability and practicality of a line
  concept, all frequencies should be multiples of a (given)  \emph{system frequency} $i\ge 2$. This can be modeled by the additional constraints
\begin{align}
	f_l=i\cdot \alpha_l &~ \forall l\in \Pool\\\
	\alpha_l \in \NN_0 &~ \forall l \in \Pool .
\end{align}

An alternative approach to model frequencies is to use indicator variables $z_l^{\phi}$  that take the value $1$ if line $l$ has frequency $\phi$, see, e.g., \cite{borndkarbstein2012,martinez2014frequency,vanderHurk2016shuttle,gatt_et_al:OASIcs.ATMOS.2022.5}.
Let $\Phi$ be the set of all allowed frequencies. We then can replace $f_l$ by $\sum_{\phi\in \Phi} \phi \cdot z_l^{\phi}$. Constraints 
\begin{align}
\sum_{\phi \in \Phi} z_l^{\phi}\le 1 &\forall   l \in \Pool
\end{align}
ensure that at most one frequency $\phi \in \Phi$ is assigned to each line $l$.

There are several lines of reasoning behind the choice of the set $\Phi$.
When the line plan to be constructed is supposed to serve as a basis for a regular timetable, i.e., a timetable in which the headways between two consecutive
trips of the same line are exactly $h_l=T/f_l$, one can argue (compare, e.g., \cite{GKIOTSALITIS2022103492}) that candidate frequencies should be divisors of $T$, so that headways and thus scheduled time points are integer. 
On the other hand, only lines with frequencies that are powers of the same base can be scheduled effectively to minimize transfer times while preserving regularity.  E.g., for a planning period of $T=60$ minutes, transfers between two lines $l$ and $l'$ with frequencies $f_l=2$ and $f_{l'}=4$, respectively, can be easily scheduled such that short transfers are possible twice an hour. For line $l$ with $f_l=2$ and $\tilde{l}$ with $f_{\tilde{l}}=3$, on the other hand, we can achieve at most one short transfer if we want to preserve regularity of the schedule. 

\cite{BornHK13} propose to strengthen line planning models 
by introducing the concept of \emph{frequency configurations}, i.e., by enumerating which combinations of line frequencies can be chosen to fulfill constraints \eqref{constraint-costmodel}. They show that adding such inequalities to transit line planning models can improve their integer linear formulation and speed up computation times.

\subsection{Transit line planning performance indicators}\label{sec-performance-indicators}
A multitude of demand- and supply-centered performance indicators and combinations thereof are used to evaluate the quality of line concepts and line plans, and (to a lesser extent) also line pools in the literature.  
All indicators can be used as a constraint (both \emph{budget} and \emph{service level} constraints are common) or as (part of) the objective function of an optimization problem.

\subsubsection{Supply-centered indicators}
 
A simple and common way to consider cost in transit line planning  is to introduce a cost parameter $\cost_l$ that models the operating cost for one trip of line $l$ and is often assumed to be proportional to the length of line $l$. Using frequency variables $f_l$ for a line pool $\Pool$, the total cost of a line concept can then be computed as $\sum_{l\in \Pool} \cost_l f_l$, as it is done in \eqref{costmodel}.
These cost are referred to as \emph{variable} cost, while $\costvariable_l$ that are
related to the introduction of a new line are referred to as \emph{fixed} cost.
For modeling fixed cost we need the indicator variables $y_l$ and receive
$\sum_{l\in \Pool} \costvariable_ly_l$.
When different types of vehicles are considered, these cost may depend on vehicle type as well, compare, e.g., \cite{gatt_et_al:OASIcs.ATMOS.2022.5}.

In models that do not use a given line pool but construct lines as part of the optimization step, fixed or variable line cost can be modeled by adding up $\textup{cost}_a$ of all PTN arcs $\kante$
contained in line $l$, i.e., $c_l:=\sum_{a \in L}\textup{cost}_a$ compare, e.g., \citep{guan2006simultaneous,BGP07}.

Energy cost are for a fixed vehicle type in a non-mountain region proportional to kilometers driven \cite{sadrani2022optimization}.
Emissions also play an increasing role in transit planning, see, e.g, \cite{bel2018evaluation,cheng2018minimizing,duran2020considering}. They also depend on the vehicle type and on the distance covered and can thus be modeled in the same way as (energy) costs, compare, e.g., \cite{beltran2009transit,pternea2015sustainable,jha2019multi,duran2019transit,duran2020considering}. \cite{pternea2015sustainable} remark that the use of electric vehicles may require additional infrastructure such as charging stations and include the cost of locating these into the objective function. 
\smallskip

In absence of more intricate cost models, variable costs are also used to estimate vehicle acquisition or maintenance cost, since these roughly correlate with the distance driven.
In order to correctly model these cost types, an estimation of the number of vehicles
to operate the line concept is needed. This is normally done under the assumption of
line-pure vehicle schedules.
In this case, the number of vehicles needed can be estimated as  $\sum_{l \in \Pool} \frac{{\rm time}_l \cdot f_l}{T}$ , where  ${\rm time}_l$ denotes the time needed to complete a trip and $T$ is the period.

\cite{PSSS17} argue that the numbers of vehicles needed should already be considered during the construction of a line pool: If the distance covered by a line is such that the line needs (slightly less) than a multiple of a planning period length to complete a trip forth and back, it can immediately be employed after finalization one trip (in a line-pure vehicle schedule), leading to higher utilization of the available rolling stock.
In the context of a subline frequency setting problem \cite{GKIOTSALITIS2022103492} explicitly determine the number of vehicles for every subline based on the length of the sublines and their frequencies. The number of vehicles needed can then be added to
the objective function.
Also \cite{canca2019integrated} and \cite{de2017railway} (in a railway context)
compute the cost for operating the vehicles based on estimated vehicle numbers.
\smallskip

All above-described costs (and emissions) are additive at the line level, i.e., each line or line-frequency combination is attributed individually with a cost and/or an emission value. The overall cost (or emissions) of a line concept are then computed as sum over the cost per line.

\subsubsection{Demand-centered indicators}

Indicators for connectivity and directness from network science can be used
to analyze line plans in transit, compare, e.g., 
\cite{luo2019integrating,fielbaum2020beyond}.
However, optimization models in transit line planning use more specific
indicators which are described in the following.

Passengers' \emph{travel time} is one of the key service quality measures used in
transit line planning, and is also used
in order to distribute passengers in the network (compare Section~\ref{sec-traveltime}). 
We distinguish between the following \emph{components of travel time}:
\begin{itemize}
 \item \emph{Driving time} is the time spent in a driving vehicle.
  Usually it is computed based on the kilometers driven and on the speed of the vehicle.
\item We refer to the time that a passenger spends in a standing vehicle at a
  stop where she does not board or alight as \emph{waiting time}. 
\item \emph{In-vehicle time} aggregates driving time and waiting time.
\item The \emph{transfer time} is the time between arriving at a station in
  one vehicle, and departing from it in a different one. In contrast to driving
  and waiting time, the transfer time depends crucially on the timetable.
\item \emph{Adaption time} is the time passing between the moment that a
  passenger desires to depart, and the time at which she actually departs.
Adaption time is particularly relevant in transit systems that operate with a high frequency where passengers do not consult a timetable prior to departure. 
  Including adaption time into the travel time can, however, also be seen as a means to penalize infrequent connections between highly demanded OD-pairs, as they lead to high adaption times  \cite{hartleb2023modeling}.
\item \emph{Access time} is the time for walking to the first or from the last
  station. It can be computed if the demand is given in districts or zones also
  considering the choice of the first and last station into the route choice
  decisions, compare \cite{dell2006bi,shimamoto2012optimisation,szeto2014transit,klier2015urban}.
  It is often neglected in transit line planning.
\end{itemize}
\smallskip

Many line planning papers use travel time or its components as objective function or as
service level constraint. Some models minimize a \emph{weighted sum} of the components, as it has been shown that the value of different time components is different in the passengers' perception, see, e.g., \cite{wardman2004public}.
When the sum over all passenger travel times is used as an objective, this is normally combined with requiring that all passengers must be transported, 
because not transporting a passenger would lead to a reduction in travel time.

Besides travel time, \emph{coverage} is a second important demand-related performance indicator in transit line planning. The definition of coverage boils down to the question of how many passengers of all OD-pairs are going to travel by transit. 

Many transit line planning models implicitly assume that a \emph{system split} (see, e.g., \cite{Oltrogge94}) is precomputed. That is, the OD passenger numbers that are taken as input do not represent the full amount of people traveling from origin to destination, but the fraction of them that wants to travel by public transport.
Based on such a system split, many line planning models require full coverage,
i.e., that \emph{all} passengers of
all OD-pairs are transported. When travel time is not considered in the objective function, such models require that the routes provided for the passengers do not deviate too much from shortest routes, i.e., a service level constraint with respect to travel time is introduced to avoid that passengers with large detours, who would in reality probably not take transit, are counted as covered. 

\emph{Direct traveler models} consider a passenger covered if there is a line that connects her origin to her destination without a large detour, see, e.g. \cite{BKZ96,Bussieck98, suman2019improvement}, while in \cite{bertsimas2021data} passengers may transfer at most once.

In contrast to such approaches, in \emph{coverage by transit}
the existence of other modes (e.g. the private car) is recognized. Only if the travel option using transit is good enough compared to the other existing mode, an OD-pair is counted as covered, see, e.g., \cite{LMMO04,Cadarso17,bertsimas2021data,klier2015urban,tirachini2014multimodal,beltran2009transit,LMO05} 
This is often framed as integrating \emph{mode choice} and will be explained in further detail in Section~\ref{sec-modechoice}.

In some models, 
instead of coverage the related indicator \emph{revenue} (as a function of coverage) is used. E.g., \cite{gkiotsalitis2022optimal,hartleb2023modeling} consider revenue from ticket prices, which are proportional to the distance between the passengers' origins and destinations.   
\cite{de2017railway, canca2019integrated} additionally include a subsidy per transported
passenger into the revenue estimation. An advantage of using revenue is that it can easily be aggregated with the supply-aggregated indicator \emph{cost} to obtain \emph{profit} as an objective function.

Also \emph{comfort}-oriented indicators are occasionally used. To ensure passenger satisfaction, besides a criterion on the travel time, \cite{gatt_et_al:OASIcs.ATMOS.2022.5} require that the buses are not filled beyond a certain level. 
\cite{cats2019frequency} include deteriorating comfort in the case of overcrowding into their simulation-based evaluation of line plans and corresponding passenger routes.
\medskip

In contrast to most supply-centered indicators that can be attributed
to individual lines, most of the demand-centered indicators can be decomposed
per passenger, but 
cannot be decomposed per line. Therefore, it is difficult or impossible to apply them to assess individual lines.
A possible approach to evaluate individual lines from a demand perspective is to consider the marginal benefit that the addition of a line has on an existing line concept, compare \cite{VerEngPhiVansteen21}.

\section{Modeling passengers and capacity}\label{sec-pax-and-cap}

In the previous section we have shown how lines and their frequencies are defined
and can be modeled. Looking at the \emph{demand-centered} performance indicators, it became clear that the routes for the passengers are
crucial for defining \emph{coverage} and \emph{travel time} indicators.
Moreover, the routes of the passengers are important to evaluate if the capacity provided
 by a line is sufficient to transport all of them.
We first describe three levels of detail in which passenger routes can be modeled and
aggregated: link-based, line-based and trip-based. We then show
for each of these aggregation levels 
how the capacity and the passenger routes are connected in
Section~\ref{sec-capacity} and how travel time can be estimated
in Section~\ref{sec-traveltime}. We finally
classify the various models for determining the passenger routes in Section~\ref{sec-demand} and end with some remarks on the mode choice.

\subsection{Modeling passenger routes in three levels of detail} 
\label{sec-aggregation}

We start by describing three levels of aggregation
which are common. For all of them, let 
$\OD$ be the OD-matrix, $\ODset$ be the set of OD-pairs
with $\OD_{s,t} > 0$ for $s,t \in V$ and let $Q_l$ be the capacity of
vehicles operating on line $l \in \Pool$.

\paragraph{Passenger routes on link level.}
Here, passenger routes  $\Pptn_{st}$ are determined in the PTN, that is, a route is a sequence of links in the PTN. The link level is the least 
detailed level. It neglects information about the lines the passengers use.

\paragraph{Passenger routes on line level.}

Here, we determine the routes of the passengers in the PTN together with the lines
to be used. 
This can be efficiently done in the \emph{Change \& Go
network} CGN, originally introduced in \cite{ScSc06a}. See Figure~\ref{fig-capacity} for
a small public transport network and its corresponding Change \& Go network.

\begin{notation}
  The \emph{change-and-go network} (CGN) $\cgn=(\cgnnodeset,\cgnarcset)$ consists
  of nodes which represent station-line pairs for each line $l$ and all the stations $v$
  it passes:
  \[ \cgnnodeset=\{ (v,l): v \in V, l \in \Pool \mbox{ and } v \in l \} \]
  The arcs set consists of \emph{transfer (=change) arcs} and \emph{driving (=go) arcs}  
  \[ \cgnarcset=\{ (\kante,l):=((v,l),(u,l)): \kante=(v,u) \in A, v,u \in l\}
  \cup \{((u,l_1),(u,l_2)): u \in l_1 \cap l_2\} \]
  which represent driving activities of $(\kante,l)$ of line $l$ on arc
  $\kante \in A$ and possible
  transfer activities for passengers between lines $l_1,l_2$ in station $u$.
  For routing an OD-pair $(s,t)$, the CGN is
  extended by adding $(s,0)$ and $(t,0)$ as nodes and
  connecting $(s,0)$ to all nodes $(s,l) \in \cgnnodeset$ and $(t,0)$
  to all $(t,l) \in \cgnnodeset$.
  The arcs in the CGN are labeled by an estimate of the travel time
  for the respective driving or transfer activity that allows a shortest path computation for determining passenger routes.
\end{notation}  

Note that the CGN can only be built if a line pool is available.

A \emph{route} $\Pcgn_{s,t}$
from $(s,0)$ to $(t,0)$ in the CGN represents a path which not only
contains the sequence of stations but also the lines to be used.
There exist a number of similar network models, compare, e.g., 
\cite{nguyen1988equilibrium,CancMautUru15}. 
Some papers use line-based aggregation levels without
the CGN, e.g. \cite{NaJe08}.

\paragraph{Passenger routes on trip level.}
A line $l$ with a frequency of let's say $f_l=4$ consists of four \emph{trips}
per period. Let ${\cal T}_l$ be the trips belonging to line $l \in \Pool$. To estimate travel time and capacity even more accurately than on the line level, we want
to specify not only the routes but also the trips which
are used on the passengers' journeys.
This is possible if line planning is integrated with scheduling decisions. In this case,
passengers can be assigned to the specific \emph{trip} to be used and not only
to the line. The exact trips a passenger uses can be computed as a route $\Pean_{s,t}$
in the so-called \emph{event-activity network (EAN)} which contains events
for each arrival or departure of every trip at every station, together with their
accurate arrival and departure times. Transit line models including timetabling and
passenger routing on trip level exist, e.g., the single
corridor frequency assignment model by
\cite{sadrani2022optimization}, or the model integrating line planning,
timetabling, and passenger routing by \cite{philinediss,SchiSchoe21}, but they
are computationally hardly tractable.

\subsection{Modeling capacity } \label{sec-capacity}

While line generation models do usually not consider capacity considerations, 
providing sufficient capacity is normally considered in line planning and
in frequency setting.
The two decisions that impact the overall capacity of line $l \in \Pool$
per period are the capacity $Q_l$ of vehicles serving the line and its frequency $f_l$.
In many models, the vehicle capacity $Q_l$ is considered to be given and identical for all lines. 
Some models allow to choose among different line types with vehicle
capacity being one distinguishing property of different lines
(besides, e.g., speed), or different vehicle sizes can be associated to
different lines, see, e.g.,
\cite{dell2012optimizing,cats2019frequency,gatt_et_al:OASIcs.ATMOS.2022.5}.

The main difference between models for transit line planning
comes from the level of capacity aggregation used, which should be
aligned with the level of detail used for the passenger routes. 
The joint determination of line \emph{frequencies} and vehicle \emph{capacities} is common in frequency setting, as both, increasing frequencies and increasing capacities,
are means to adjust line capacity to passenger demand. While in link-based and
line-based approaches, all vehicles of the same line have the same capacity
\cite{dell2012optimizing,cats2019frequency,gatt_et_al:OASIcs.ATMOS.2022.5}
it may make sense in trip-based models 
to assign vehicles of different sizes to the same line, to cater for fluctuations in demand, compare \cite{sadrani2022optimization}.
In the following we formulate the relation between frequency, capacity and demand for the three aggregation levels link-based, line-based and trip-based.

\paragraph{Link-based capacity aggregation.}
For the link-based aggregation level we use the links of the PTN 
defined in Notation~\ref{nota-PTN}. In this level of detail 
we can aggregate the frequencies of all lines for each specific link and require that
the capacity on the link
is sufficient for the passengers that use this link,
independent of the specific lines they take. 
Constraint
\begin{equation}
  \label{demandPTN}
  \sum_{l \in \Pool: a \in l} f_l \cdot Q_l \geq \sum_{(s,t) \in \ODset: \kante \in \Pptn_{s,t}} \OD_{s,t}=d_a
  \mbox{ for all } \kante \in A
\end{equation}  
makes sure that the demand is covered on every link $\kante \in A$ of the PTN when
solving the line planning problem. The cost model \eqref{costmodel} is an example
which uses the link-based aggregation.

\paragraph{Line-based capacity aggregation.}
The link-based aggregation neglects which lines passengers use.
Even if the capacity on a link $\kante \in A$  is large enough to transport the demand $d_a$
it might be possible that not all passengers can use their preferred lines. 
To ensure that the capacity \emph{per line} is sufficient to transport all passengers
we use routes $\Pcgn_{s,t}$ in the CGN and receive 
\begin{equation}
  \label{demandCGN}
  f_l \cdot Q_l \geq \sum_{(s,t) \in \ODset: (\kante,l) \in \Pcgn_{s,t}} \OD_{s,t}=:d_{\kante,l}
  \mbox{ for all } \kante \in A, l \in \Pool
\end{equation}  
making sure that the number of passengers using line $l$ on arc $\kante \in A$ 
does not exceed the \emph{line capacity} $f_l\cdot Q_l$ on every arc $\kante \in A$
of the line, compare,
e.g., \cite{guan2006simultaneous,ScSc06a,gatt2021column}.
In other words, it is ensured that over the whole period $T$, capacity
on each line $l$ is sufficient to transport the passengers
assigned to the line. \eqref{demandCGN} also
confirms that no demand can use a line $l \in \Pool$ which is not
selected,  i.e., $l \not\in \cL$. 

\paragraph{Trip-based capacity aggregation.}
When passenger routes are detailed on trip level, capacity can be modeled on a trip
level as well. That is, for each arc $\kante \in A$ 
and each trip operated on $\kante$ 
we count the passengers that are assigned to this trip and bound their number by
the number of places in the vehicle. With ${\cal T}$ being the set of all trips,
we receive
\begin{equation}
  \label{demandEAN}
  Q_l \geq \sum_{(s,t) \in \ODset: (\kante,l,\textup{trip}) \in \Pean_{s,t}} \OD_{s,t}=:d_{\kante,l,\textup{trip}}
  \mbox{ for all } l \in \Pool, \textup{trip} \in {\cal T}_l, \kante \in A, 
\end{equation}

The following example illustrates why trip-based capacity modeling may be desirable.

\paragraph{Illustration.}

\begin{figure}
\begin{center}
  \includegraphics[width=0.42 \textwidth]{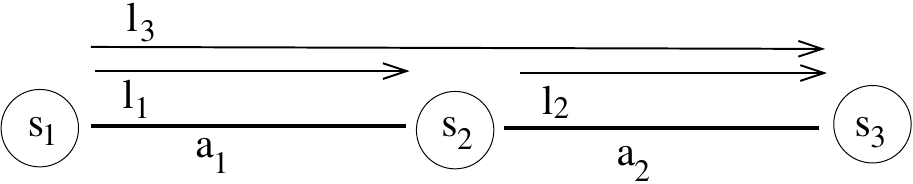}\hspace{1cm}
  \includegraphics[width=0.42 \textwidth]{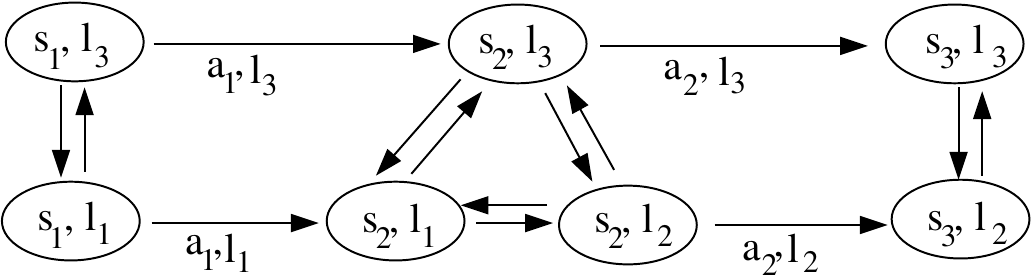}
\end{center}
\caption{The public transport network and the Change \& Go network for the illustration.}
  \label{fig-capacity}
\end{figure}

We illustrate and compare the three aggregation levels in a small example. Consider
the network with three stations, $s_1,s_2,$ and $s_3$ and three lines
$\cL=\{l_1,l_2,l_3\}$
depicted in Figure~\ref{fig-capacity}.
Assume that 120 persons wish to travel from $s_1$ to $s_3$.
We want to evaluate a solution in which the
lines $l_1$ and $l_3$ operate with a frequency of $f_{l_1}=f_{l_3}=1$
and capacities of $Q_{l_1}=Q_{l_3}=60$ while line $l_2$ runs with the double
frequency of $f_{l_2}=2$, but only uses smaller buses with half of the
capacity, i.e., $Q_{l_2}=30$.

\begin{itemize}
\item In the link-based aggregation level we determine the traffic loads of
  both links $\kante_1=(s_1,s_2), \kante_2=(s_2,s_3) \in A$ as
  $d_{\kante_1}=d_{\kante_2}=120$. We receive
  \begin{eqnarray*}
    \sum_{l \in \Pool: \kante_1 \in l} f_l Q_l & = & 60 + 60 \geq 120=d_{\kante_1}, \\
    \sum_{l \in \Pool: \kante_2 \in l} f_l Q_l & = & 60 + 30 + 30 \geq 120=d_{\kante_2}, 
  \end{eqnarray*}
  so on a link-based aggregation the solution is feasible with respect to \eqref{demandPTN}
  although we may expect
  overcrowded buses on line $l_3$ since most passengers would prefer a direct
  connection over a connection with a transfer.
\item The line-based aggregation would either route all passengers along line $l_3$
  yielding an infeasible solution since capacity of $l_3$ is only limited to 60 persons.
  With a capacity-aware routing model, 60 persons would use the path
  $\Pcgn_1=((s_1,l_3),(s_2,l_3),(s_3,l_3))$ and 60 persons would be assigned to path
  $\Pcgn_2=((s_1,l_1),(s_2,l_1),(s_2,l_2),(s_3,l_2))$. With these two paths,
  the capacity constraints \eqref{demandCGN} are satisfied, since
  \begin{eqnarray*}
    f_{l_1} Q_{l_1} = 60 \geq 60= d_{\kante_1,l_1},\\
    f_{l_2} Q_{l_2} = 60 \geq 60= d_{\kante_2,l_2},\\
    f_{l_3} Q_{l_3} = 60 \geq 60= \max\{ d_{\kante_1,l_3}, d_{\kante_2,l_3}\}
  \end{eqnarray*}
  and the solution would be classified as feasible. 
  However, the 60 passengers arriving at station $s_2$ will all try to use
  the same bus of line $l_2$ (namely the first one that arrives) leading to an
  overcrowded bus --- while the next bus of line $l_2$  would be not used at all.
\item The trip-based model assigning passengers to the individual trips is the only of these three which
  reflects capacity correctly in this situation: 
  For the trips of lines $l_1$ and $l_3$ we see that their capacities are sufficient to transport all
  passengers assigned to them in one trip.
  However, there are two trips (let's say $\textup{trip}_1$ and $\textup{trip}_2$)
  belonging to line $l_2$. All 60 passengers arriving at
  station $s_2$ with line $l_1$ will board the first of these two trips,
  say $\textup{trip}_1$. Then we have
  \[ Q_{l_2}=30 < 60=d_{\kante_2,l_2,\textup{trip}_1}, \]  
  i.e., \eqref{demandEAN} is not satisfied. 
\end{itemize}

\subsection{Estimating passengers' travel time}
\label{sec-traveltime}

Estimating travel time is key to most demand distribution models (compare Section~\ref{sec-demand}): the travel time of a certain route determines how many passengers take it, and the amount of passengers on a route in turn determines whether capacity constraints are satisfied (compare Section~\ref{sec-capacity}). 
Travel time and its components are also relevant performance indicators to evaluate a line concept, and, as such, have been described in Section~\ref{sec-performance-indicators}. 
This section examines how travel time can be estimated based on the level of detail
used.

\paragraph{Estimating travel time in case of trip-based aggregation.} 
An exact estimate of planned travel time can be given in presence of a timetable.
If there are no delays, the departure time of the first boarded vehicle at the origin and the arrival time of the last vehicle of a route at the destination can be read from the timetable. Sometimes, even denied boarding is modeled and accounted for in the transfer and adaption time, compare  \cite{sadrani2022optimization}.  
If the desired departure time of a passenger was known (normally, this is not the case in line planning), also adaption time could be computed exactly. Otherwise, it must be estimated, e.g., as in \eqref{transfertime}.

\paragraph{Estimating travel time in case of line-based aggregation.}
When passenger routes are modeled line-based, \emph{driving time} can be estimated based on line speed and distance and assigned to the arcs of the CGN. 
In a few models for headway-based transportation, road congestion is explicitly
modeled. E.g., \cite{tirachini2014multimodal} allow passengers to choose between
different modes, where both, bus and car traffic, contribute to crowded streets.

Sometimes, \emph{waiting time} is included in the in-vehicle time as a constant supplement to the driving times,
 \cite{szeto2014transit} assume that waiting time depends
linearly on the number of passengers already in the vehicle (modeled as passengers
assigned to the line on that link divided by the frequency),
 but usually the vehicles' stopping times are considered as marginal and hence not included in the travel time. 

 In absence of a timetable, transfer times cannot be predicted exactly, but need  to be estimated.
 To assign transfer time as label to the CGN, 
the transfer time can be approximated either  by using a fixed  penalty term \cite{guan2006simultaneous,ScSc06a}  or by using the frequency of the departing line.
An approximation for the transfer time in the line-based model
is the average expected transfer time
\begin{equation}
 \label{transfertime}
\frac{T}{2 \cdot f_l}.
\end{equation}  
The same formula can be used for the  \emph{adaption time}.
When using binary indicator variables, 
this is equivalent to 
$\sum_{\phi \in \Phi} \frac{T}{2\phi}z_l^{\phi}$,  compare, e.g.,
\cite{bull2015optimization,sadrani2022optimization}). 

\cite{gatt_et_al:OASIcs.ATMOS.2022.5} follow the same approach when the
frequency is high (more than 6 per hour), but argue, that for low frequencies
passengers will mind the timetable to avoid adaption time
and transfers will be synchronized. Therefore, for lower frequencies
they assume a transfer and adaption time of $5$ minutes.
Interestingly, a similar result also holds in the absence of regularity constraints:
\cite{GKIOTSALITIS2022103492} prove that in the case of $k$ periodically
scheduled but un-synchronized (i.e.: irregular) vehicle departures the expected
adaption/waiting time to board one of them is $\frac{T}{k+1}$, which
even holds for the case of irregular departures.  

\paragraph{Estimating travel time in case of link-based aggregation.}
Routes that are modeled on link-based aggregation level are paths in the PTN and hence
do not encode which lines are taken.
Most models (see, e.g.,\cite{BGP07,borndkarbstein2012,gatt2021column,bertsimas2021data}) assume constant line speed so that driving times do not depend on the chosen lines and can hence be estimated  without knowing the lines used by the passengers. Transfer times are mostly neglected in the link-based aggregation, see, e.g., \cite{BGP07}.
A possibility is to distinguish between direct routes and routes with at least one transfer and penalize the latter ones \cite{borndkarbstein2012,gatt2021column,bertsimas2021data}. In some direct traveler models, routes in the PTN are associated with the fastest line (or just any line) that covers the route, and a passengers' travel time is given between the passenger's origin and destination, compare, e.g., \cite{suman2019improvement}). In link-based aggregation, adaption time is normally not included.

\subsection{Modeling passengers' route choice} \label{sec-demand}
\begin{figure}
\begin{center}
\includegraphics[width=1 \textwidth]{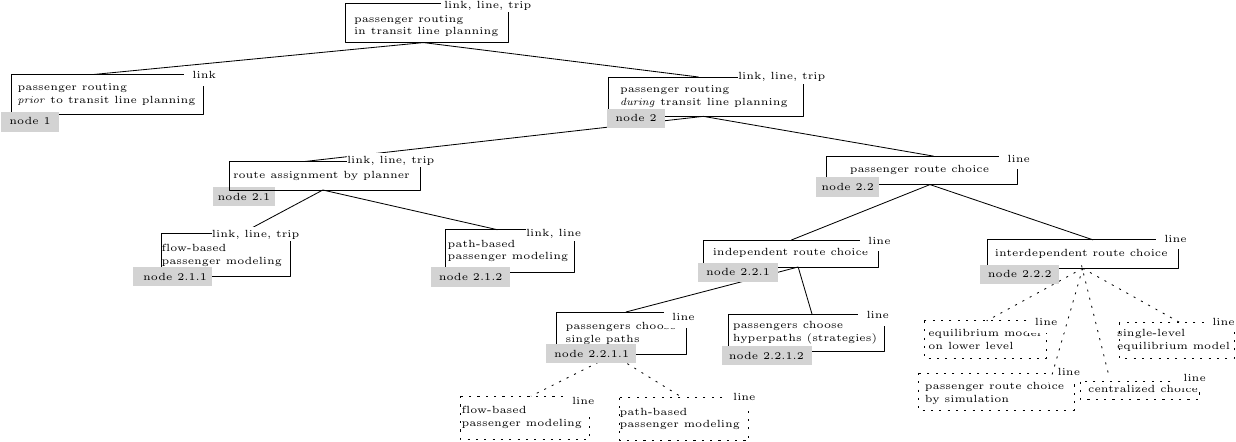}
\end{center}
\caption{Models for passengers' route choice and their possible aggregation levels.}
  \label{fig-baum}
\end{figure}

To determine travel time or to evaluate the capacity constraints
\eqref{demandPTN}, \eqref{demandCGN}, or \eqref{demandEAN} in the different
aggregation levels, passenger routes need to be determined. 
To describe the many route choice models to do this, we use the structure shown in the
decision tree in Figure~\ref{fig-baum}.

We first distinguish between approaches that distribute passengers over the network
\emph{prior} to the optimization step
and models that aim to integrate the route choice of the passengers into the optimization step.

\paragraph{Passenger routing \emph{prior} to transit line planning (node 1)}
In order to distribute OD-demand prior to the construction of a line concept, shortest-path algorithms or other assignment models, see, e.g., \cite{de2011modelling} on the PTN can be used.
The predetermined passenger distribution can be used to compute 
\emph{traffic loads} on each link of the PTN, see \eqref{trafficload}. This is
a common approach, e.g., used in \cite{BuLiLue04,CDZ98,GHK06}.
Other models use precomputed PTN-paths for the OD-pairs that 
need to be enabled in order to travel, see \cite{BKZ96}.

\paragraph{Passenger routing \emph{during} transit line planning (node 2)}
We now turn to models that do not consider the passenger routes as given, but acknowledge that they depend on the transit line planning decisions to be made in the optimization.
We distinguish between optimization models where passengers are \emph{assigned} to routes by a central decision maker (node 2.1) 
and models that assume that passengers choose a route by solving their own
optimization problems (e.g., a shortest route with respect to travel time). 

\paragraph{Route assignment assignment by planner (node 2.1)}
We first review models that aim at a \emph{system optimum}, i.e., the sum of travel times (or any other performance indicator) is to be minimized by a central decision maker who makes both decisions, on the supply (lines and
frequencies) and on the demand (passengers routes). The decision maker 
optimizes \emph{her} objective function which is not necessarily fully aligned with the passengers individual
objective functions. I.e., such models do not forbid to assign individual passengers to routes that are
suboptimal from the perspective of the individual passenger, if this benefits the overall objective function.
Since in this case, transit line planning reduces to a  single level problem,
it can be often approached by integer-programming.

Other publications explicitly choose passenger routes during the optimization process and will be reviewed below. For doing so, some rely on a precomputed \emph{candidate set} of routes for each OD-pair, while others use \emph{flow constraints}.  

\paragraph{Flow-based passenger modeling (node 2.1.1)}
\emph{Flow-based} models describe the passenger routes as multi-commodity flows between origins and destinations. This approach can be applied on all aggregation levels for passenger routes (compare Section~\ref{sec-aggregation}). That is, in the description that follows, the term \emph{network} can refer to the PTN (e.g., in \cite{FHSS17,LMMO04,Cadarso17}), CGN (see, e.g.,  \cite{NaJe08,ScSc06a}) or EAN (see, e.g., \cite{SchiSchoe21,LiTang23}). 

Such models use flow variables $x_{\kante}^{\odpair}$ that indicate how many passengers of an OD-pair $\odpair$ travel on an arc $\kante$ of the network. Flow constraints are modeled as
\begin{align}
\theta x^{(\org,\dest)}=\paxnumber^{\odpair} \ \ \forall \odpair \in \ODset \label{eq-ptnflow}
\end{align}
 where $\theta$ represents the node-arc incidence matrix of the network and $\paxnumber^{(\org,\dest)}$ is the number of passengers on OD-pair $\odpair$.
Models that use passenger routing with flow constraints therefore minimize total (or, equivalently, average)
travel time. Note that
passengers may be assigned to routes with long travel time if this benefits the objective function of the model.
In order to avoid this, one can impose bounds on the maximal travel time for each OD-pair, depending on the geographical distance of origin and destination or on the length of a shortest route, see, e.g., \cite{FHSS17}.

\paragraph{Path-based passenger modeling (node 2.1.2)}
As alternative to flow-based models, some models \emph{precompute choice sets} of routes on link or line level for each OD-pair $\odpair$, and assign the passengers of $\odpair$ to one or several of these routes during the optimization step. In most models, only routes of a certain quality, e.g., that do not exceed the length of a shortest path by too much, are included in the choice set. From a modeling perspective, this approach has the advantage that passengers cannot be assigned to 'bad' routes for the sake of the social optimum, even if the objective function
minimizes not travel time, but, e.g., costs.
On the other hand, route choice is more constrained in these models, in particular if the choice set is small.
To overcome this problem,  \cite{BGP07} use a column generation approach to iteratively \emph{generate} suitable
OD-paths in the PTN.

A challenge is that in particular in dense line pools, the number of routes with good
travel time for an OD-pair $\odpair$ can be very big. In the CGN this may happen
even if there is only one 
path from $\org$ to $\dest$ in the PTN.
Therefore, preprocessing the choice set is a promising approach. E.g., \cite{hartleb2023modeling}
restrict to routes with only one transfer and set tight bounds on the travel time.

On the contrary, \cite{klier2015urban} consider an even bigger choice set 
where each CGN path corresponds to multiple elements - one for each possible frequency combination for the lines used on the CGN path. This allows to precompute travel time including frequency-dependent adaption and transfer times.

\paragraph{Passenger route \emph{choice} (node 2.2)}
We now review approaches that explicitly model route \emph{choice} of passengers.  
Such models are often represented as bilevel optimization models with transit line planning decisions (route generation, selection, and frequency setting) on the upper level and passenger routing on the lower level. 

\emph{Solution approaches} for these kind of models can be roughly categorized as follows: The first class of approaches computes transit line planning decisions and routing simultaneously, e.g., by solving a mixed-integer program that encodes both transit network design decisions and decisions on passenger routes in variables of the same problem. 
The other class of approaches follows the hierarchy suggested by the bilevel representation and takes transit line planning decisions first, and routing decisions second. In many heuristics and metaheuristics this process is iterated, interpreting the computation of a routing for a given line concept  as an evaluation mechanism and using it to guide the solution process.
The first approach has the advantage that (if the resulting models can be solved to optimality) a global optimum is found, while the second approach can get stuck in local optima. On the other hand, the routing models of the first
approach need to be simple enough to be used
within a mixed-integer linear program, while the stepwise treatment of line planning decisions and routing in the second approach has the advantage, that line planning decisions are already fixed when a routing is computed, thus more sophisticated routing models can be used.
\smallskip

We distinguish further by the type of problem considered on the lower level.

\paragraph{\emph{Independent} route choice (node 2.2.1)}
We first turn our attention to models where the routing problem of each passenger can be formulated as an optimization problem based on the upper level decisions, but independent of the choices of other passengers.
Here, we distinguish between models where the solutions to the lower level of optimization problems for an individual
OD-pair are \emph{routes}, and models where lower level solutions are \emph{strategies}.

\paragraph{Passengers choose \emph{single paths} (node 2.2.1.1)}    
These models are motivated by the consideration that, as soon as the decisions on line routes and frequencies (and possibly a timetable) are taken, transfer times and adaption times can be computed as described in Section~\ref{sec-traveltime}, and finding a travel-time minimal line-based (trip-based) route 
from an origin to a destination corresponds to solving a shortest path problem in the CGN (or the EAN). 
 
If no capacity constraints are considered, and if the upper level objective is the sum of the individual OD-pair objective functions of the lower level optimization problems, in a social optimum all passengers are assigned to individually optimal paths, and thus, the bilevel problem can be transformed to a single level problem with route assignment by planner as described in node 2.1, see \cite{mariediss}. 
However, in general this does not hold. 
\cite{mariediss,goerigk2017line} remark that in the presence of capacity constraints, 
the passenger routes estimated by flow constraints 
may not follow shortest paths. They argue that this may lead to capacity violations in reality, when passengers do behave rational and choose shortest routes instead of following the routes predicted by the model. 
Therefore, they propose to model transit line planning as a bilevel program, where the operator's objective constitutes the upper level and passengers' route choice (using line-based aggregation) constitute the lower level. By introducing additional \emph{shortest-route-indicator} variables, or by dualizing the lower level problem, the bilevel problem can be transferred to a mixed integer program, which is, however, hard to solve
is followed in the route-selection based frequency setting model by \cite{gatt_et_al:OASIcs.ATMOS.2022.5} where the authors use an optimal value function to assign passengers to the shortest
path that is enabled by the line concept. 

\paragraph{Passengers choose \emph{hyperpaths (strategies)} (node 2.2.1.2)}  
The above-described passenger routing models 
have in common that they use (adaption and) transfer time estimates (either frequency-based, fixed penalties, or no transfer time) in order to compute travel-time minimal routes and then assign passengers to these routes. Unfortunately, the travel time of a route may deviate considerably from the estimates made.

For this reason, some authors argue that assignment of OD-pairs to routes is not the best way to estimate passenger routes choice in line planning. 
The approach that they propose is motivated by a routing paradigm for headway-based 
transportation, described first as the \emph{common lines problem} by \cite{chriqui1975common}:  When several lines could be boarded to go from origin to destination (possibly with transfers), in absence of a timetable
a passenger would not decide prior to departure which line to take, but make her decision depend on, e.g., which line arrives first. 
That is, instead of choosing a route prior to departure, passengers
choose a \emph{strategy} based on the line plan, i.e., a set of rules
that the user will apply in an online fashion to travel through the
system, see \cite{spiess1989optimal}. 

This strategy can be represented as an acyclic directed graph or hyperpath in a suitably chosen network  (compare, e.g., \cite{constantin1995optimizing,shimamoto2012optimisation,CancMautUru15}), where each node with more than one outgoing arc represents a decision on how a passenger will continue her route. 
Based on line frequencies, probabilities for each OD-path in the network
can be determined 
which allow to compute the expected travel time of
a strategy by solving a linear program, see \cite{spiess1989optimal}. 
This route choice approach has been integrated as lower level problem in a bilevel frequency setting approach \cite{constantin1995optimizing}, which has been transformed into a mixed-integer program in \cite{martinez2014frequency}.

\cite{shimamoto2012optimisation} use the computation of optimal strategies for the passengers as fitness function in a metaheuristic for the TNDFS.

\paragraph{Interdependent route choice (node 2.2.2)}
Lastly, there are models where the route choice of individual passengers does not only depend on the line plan, but is also impacted by the route choice of other passengers.

\cite{CancMautUru15} use a bilevel model where on the lower level they optimize all the passenger routes
based on strategies simultaneously.
There is also a number of publications that formulate transit line planning with equilibrium constraints on the lower level. To give two examples,
\cite{dell2012optimizing} 
formulate a joint frequency and bus size assignment problem as bilevel optimization problem. They apply a metaheuristic where, on the lower level, they use a transit assignment model proposed by  \cite{de1993transit}, which leads to \emph{equilibrium constraints} for passenger flow modeling.
A similar model is employed in \cite{dell2006bi} 
for joint stop location and frequency setting.
Within a genetic algorithm, \cite{beltran2009transit} 
use an equilibrium model where passengers can choose among three ways of traveling; by car, by public transport in conventional vehicles. 

\cite{SchieweSchieweSchmidt19} propose a model for transit line planning with the full pool where passengers choose their (line-based) routes based on a weighted sum of travel time and \emph{distributed operational cost}, i.e., when a passenger chooses a previously unused line or frequencies are increased, this is more costly than choosing a route that uses already established lines and their frequencies. They sequentially apply this routing strategy in a best-response approach to find a line concept. A similar approach is applied in \cite{gioialine} for line generation.

\cite{cats2019frequency} use \emph{agent-based simulation} to model lower-level route choice in a frequency assignment problem. 

\subsection{Mode choice}\label{sec-modechoice}
While many line planning models assume that the \emph{system split} is precomputed and that passenger demand is inelastic in the sense that the OD-matrix is given and fixed and represents all passengers who travel by transit, there are models that recognize that passengers are sensitive to the quality of transit and aim to include \emph{mode choice} in transit line planning.

Some papers propose an all-or-nothing assignment: If and only if the offered connection from origin to destination is better than an alternative mode, passengers
choose it. This includes the models which maximize coverage by transit as performance
indicator: Here, a passenger uses transit if the travel time is smaller than a (maybe
fictive) travel time in an alternative mode, e.g., the private car, see, 
\cite{LMMO04,LaporteMesaPerea10,CadMar12,CadMar16,Cadarso17}.
A link-wise all-or-nothing assignment is used in \cite{gkiotsalitis2022optimal} 
studying frequency assignment under pandemic-imposed distancing measures.
They allow to split passenger flows into passengers traveling on the line
segment, and passengers bridging the distance covered by the line segment by
other means. 

In order to estimate the number of passengers using public transit, also
\emph{logit models} based on travel time have been used. 
\cite{de2017railway} use a piecewise linearization of the logit function.
Also \cite{canca2019integrated} use a piecewise linearization of the logit
function for route choice
within a very detailed modeled transit line planning problem.
\cite{tirachini2014multimodal} formulate mode choice for three alternative modes (bus, car, and walking) using the logit model, for a frequency setting that integrates not only decision on bus capacity, but also on fare, congestion toll, fare collection system, and bus boarding policy on a linear network. They do not linearize their formulation but solve it directly using sequential quadratic programming.
\cite{beltran2009transit} consider three options to provide transport along a link: public transport with regular vehicles, public transport with green vehicles and transport by car, and employ an equilibrium model for route and mode choice. 
\cite{LMO05} use a logit model for locating a single line and solve it with a greedy
approach.
\smallskip

Others argue that besides travel time, the number of connections per period
plays a relevant role in attracting passengers to travel by public transit.
\cite{hartleb2023modeling} argue that, in a reasonable line plan, in-vehicle time and transfer time on routes between OD-pairs is negligible, and thus, the utility of using public transit depends only on the number of \emph{reasonable} routes offered between origin and destination per hour.
Also, in their direct traveler model, \cite{bertsimas2021data}
assume that the percentage of passengers of an OD-pair choosing
public transit depends on the frequency of the connecting line offered.
Assuming that public transit demand  depends only on the quality of the best available \emph{route bundle}, \cite{klier2015urban} are able to \emph{precompute} the OD-demand.

\section{Transit line planning under uncertainty}
\label{sec-uncertainty}

The goal of uncertain optimization in transit line planning is to find a line concept which
is still good under variation of input parameters. We discuss the following sources of uncertainty.
\begin{itemize}
\item Uncertainty of the demand, i.e., how many passengers wish to travel (see Section~\ref{sec-demand-uncertainty}),
\item unforeseen link failures in the PTN (see Section~\ref{sec-network-uncertainty}), or
\item uncertainty of the driving times which may vary due to congestion, weather conditions, or
  unforeseen delays (see Section~\ref{sec-traveltime-uncertainty}).
\end{itemize}
Each realization of an uncertain parameter is called a \emph{scenario} and the set of all scenarios is called the \emph{uncertainty set} $\cU$. Sometimes, a probability distribution on $\cU$ is known. Further aspects of uncertainty can be considered, e.g., 
\cite{KonZaro08,BesKonZar09} model competition between railway operators which
do not reveal their real incentives.

\paragraph{Optimality concepts.}  
Uncertain optimization can be treated by robust and stochastic optimization
(see \cite{RObook,BirgeLouveaux97,KouYu97}). They differ in the optimality concepts
they follow.
\smallskip

\emph{One-stage robust optimization} (e.g., \cite{RObook})
hedges against the worst case. It hence aims for a solution which is feasible
for every scenario in the uncertainty set and performs best in the worst case
over all scenarios. The concept is also called \emph{strict robustness}.
\smallskip

\emph{Two-stage robust optimization} splits decisions into
\emph{here-and-now decisions} of the first stage 
which have to be made without knowing the scenario and \emph{wait-and-see decisions} 
in a second stage. The latter can be taken after the scenario is revealed.
The goal is to find here-and-now decisions that can be recovered to a good
solution even in the worst case over all scenarios. These optimality concepts
include adjustable robustness (e.g., \cite{BeGoGuNe2003}) and recovery robustness
(e.g., \cite{LLMS09}). 
\smallskip

\emph{One-stage stochastic optimization} \cite{BirgeLouveaux97}
assume that a probability distribution on the scenario set is known.
This allows to minimize the expected value of the objective function.
\smallskip

\emph{Two-stage stochastic optimization} \cite{BirgeLouveaux97}
considers a second stage in which recourse variables
reflecting wait-and-see decisions can be chosen after the scenario has been revealed.
Several goals may be considered, among them the minimization of the expected value.

\paragraph{Uncertainty sets.}
We introduce uncertainty sets which have been used in transit
network design. For illustration, consider a problem in which two parameters,
$(l_1,l_2)$ are uncertain. These may represent uncertain driving times on two links $a_1, a_2$.
\begin{itemize}
\item In the \emph{deterministic case} the values of $l_1$ and $l_2$ are known and the
  uncertainty set contains only one point, e.g., $\cU=\{(8,15)\}$.
\item A \emph{discrete uncertainty set} consists of a finite set of scenarios. These could be
  driving times under three different traffic situations, e.g., 
  $\cU=\{ (10,20), (9,17), (8,15) \}$. Probabilities may be assigned to each of the scenarios.
\item If parameter ranges can be specified independently for each parameter, 
we obtain \emph{box-wise uncertainty}. 
If the driving time $l_1$ on link $a_1$ is in $[8,10]$ and 
$l_2 \in [15,20]$ we receive $\cU:= [8,10] \times [15,20]$.
A typical distribution on $\cU$ (reflecting measurement errors) may be the multi-normal
distribution.
\item \emph{$\Gamma$-uncertainty} introduced by \cite{BertSim04}
  takes into account that in the case of box-uncertainty it may be unlikely that all
  uncertain parameters take their extreme values simultaneously.
  In the example, it might be unlikely, that both driving times take their worst
  case at the same time. This is reflected by a budget constraint with budget $\Gamma$,
  leading, e.g., to $\cU=\{(l_1,l_2): 8 \leq l_1 \leq 10, 15 \leq l_2 \leq 20, \mbox{ and }
  l_1 + l_2 \leq 25 =:\Gamma \}$.
\end{itemize}
\medskip

Uncertain optimization generates challenging optimization problems, and many of them
may be practically relevant. Nevertheless, research on transit line planning under
uncertainties is still scarce. In the following we sketch how some of the introduced
optimality concepts have been applied to transit line planning.

\subsection{Transit line planning under demand uncertainty}
\label{sec-demand-uncertainty}

Most publications assume that transit demand, as given in the OD-matrix, is given and
fixed, or that its dependence on the supply provided is deterministic
(compare Section~\ref{sec-modechoice}).
However, since exact predictions about future demand are hard to make,
the OD-matrix should be considered as an uncertain parameter.
In the papers sketched below, the goal is to design a line concept which 
still can transport all passengers even if demand increases. Box-wise 
uncertainty $\cU_{box}=\{\OD: l_{s,t} \leq \OD_{s,t} \leq u_{s,t} \ \forall (s,t) \in \ODset\}$
for given lower and upper bounds $l_{s,t},u_{s,t}$ for all OD-pairs $(s,t)$
and $\Gamma$-based uncertainty
$\cU_\Gamma=\{\OD \in \cU_{box}: \sum_{(s,t) \in \ODset} \OD_{s,t} \leq \Gamma\}$
have been used.
\smallskip

\textbf{One-stage robust optimization.}
\cite{LeeNair21} consider a strictly robust approach to transit line planning
in which they assume $\Gamma$-uncertainty. They search a line concept and passenger routes that
minimize the travel time of the passengers including a penalty for overcrowded lines in the worst
case of all OD-matrices in $\cU_{\Gamma}$. They use a line-based aggregation and model both, route assignment
by the planner and route choice by the passengers.

\textbf{One-stage stochastic optimization.}
\cite{GKIOTSALITIS2022103492} study a single-stage stochastic version of the \emph{subline frequency setting problem}, i.e., the problem to select lines and frequencies to run on a corridor. In a deterministic version of the problem, frequencies have to be chosen high enough to allow all passengers to board the next-arriving vehicle towards their destination, however, in the stochastic case this constraint may be violated and capacity violations lead to a penalty. They solve the problem by scenario sampling from a known distribution of the demand.
\smallskip

\textbf{Two-stage robust optimization.}
Box-uncertainty has been considered in \cite{AnLo16} and in \cite{pu2021two}.
Both groups of authors minimize a sum of costs and travel time
in the worst-case over all possible scenarios in the line-based aggregation
level. They use two-stage models where 
lines and frequencies are determined in the first stage
before knowing the scenario, while passenger routes are adjusted to the
realized scenario in the second stage.

In \cite{AnLo16} two models are discussed.
In the first one, passenger routes are assigned by a planner
under hard capacity limits of the lines.
In the second model, line capacity violations are allowed, but lead to discomfort
for the passengers who choose their routes interdependently according to an
equilibrium model. In \cite{pu2021two}, passenger routes are assigned
by the planner. 
Here, additionally to the adjustment of the passenger routes, the authors allow to also adjust
the routes of the lines on a limited number of links.
Their objective function also includes the travel time and costs in the first stage.

In both papers, the problems where passenger routes are assigned
by the planner 
can be solved by using the worst case in which all OD-pairs
take their maximum values. For the model with equilibrium routing presented in
\cite{AnLo16}, the authors observe a paradox: it may happen that maximum demand does
not constitute the worst-case scenario.
\smallskip

\textbf{Two-stage stochastic optimization.}
\cite{AnLo16} propose
two-level stochastic optimization models to
transit line planning with uncertain demand with box uncertainty. They
minimize the expected cost and travel times assuming a known probability
distribution of the passengers' demand which is made tractable
by sampling approaches. 
As in their robust optimization approaches described above, the first stage
concerns the decisions on lines and their frequencies and the second
stage allows to adjust the passenger routes within a line-based level of detail
when the scenario is known.
They analyze the situation where  passengers are assigned to routes by
the planner as well
as routing passengers interdependently with equilibrium flows assuming capacity
violation penalties.
Additionally, the second stage does not only adjust the routing variables but 
also allows to add flexible service (such as dial-a-ride) on higher costs for passengers
that cannot be transported by the transit lines.

\subsection{Transit line planning under the risk of link failure}
\label{sec-network-uncertainty}

If a link in the PTN fails, some passengers might not be able to travel any more or they might
suffer detours. The goal in transit line planning is to find a line concept which keeps the number
of badly affected passengers low even if links fail.
The uncertainty set $\cU$ is discrete in this setting since there is only a finite number
of links to fail. If only the failure of single links is considered we have $\cU=A$.
\medskip

We first discuss how the consequences of link failures can be evaluated.
Link failures are a typical disruption in classic network design problems and have also been
investigated in transit line planning. While many papers evaluate networks from a pure topology
point of view (e.g., \cite{derrible2010characterizing}),
\cite{SantosLaporteMesaPerea12} were one of the first who included the demand in a line-based
evaluation. For the failure of a link $a$ they propose two models: one in which
passengers need to avoid the link, i.e., they compute new (shortest) paths with
an increased travel time. The second model assumes bridging the failed link by a bus
service which also increases the travel times. 
A full-scan of all edges is performed and the maximum of them is used as robustness index.
Also \cite{Cats16} suggests a full-scan of all links, but on a link-based level.
In case of a failure, they distinguish between passengers without any path, passengers which
need a detour and undisturbed passengers. 
\smallskip

\textbf{One-stage robust or stochastic optimization.}
If passenger routes cannot be adjusted, failure of a link will destroy
some passenger paths (unless the link would not be used at all). 
The worst case would be attained at the link with the highest traffic load.
However, in real life passengers would adjust their routes which makes one-stage
optimization concepts for link failures useless
from a practical point of view and might be the reason
that they have not been considered for link failures in transit line planning.
\smallskip

\textbf{Two-stage robust optimization.}
In a two-stage approach, the lines and their frequencies are
set in the first stage while the routes of the passengers can
be adjusted in the second stage when it is known which link has
failed.
This concept has been followed in \cite{MarinMesaPerea09}, who minimize the
passengers' traveling times under link failure in the worst case.
The authors analyze travel time in the link-based aggregation level
assigning passengers to routes by the planner.
The same setting is used in \cite{LaporteMesaPerea10} but with coverage by transit
in presence of an alternative mode 
as performance indicator, i.e., they aim
to find a line concept that maximizes the coverage of
passengers in the worst case. Since the uncertainty set with $\cU:=A$
is small, the authors can
formulate an integer linear program which contains paths for all OD-pairs and all scenarios.
They show that it is equivalent to a two-player non-cooperative zero-sum game in which
the first player is the transit line planner which can choose from different line concepts in his strategy
set while the second player is an evil adversary who selects a link to fail.
The authors prove that a saddle point of the game (if it exists) corresponds to a
(strictly) robust solution of the two-stage robust transit line planning problem.
\smallskip

\textbf{Two-stage stochastic optimization.}
Given that probabilities for the link failures are available,
\cite{LaporteMesaPerea10} maximize 
the expected coverage of the line concept. Passenger routes are determined
(by the planner)
in the second stage in the link-based aggregation level, and may be 
adjusted to the link which has failed.
Analogous to their approach for minimizing the worst case, and due to the rather
  small scenario set $\cU=A$ the resulting problem can be formulated as linear integer program.

\subsection{Transit line planning under driving time uncertainty}
\label{sec-traveltime-uncertainty}

Uncertainty in the driving times may, e.g., be due to varying traffic, or
due to unforeseen delays. Hedging against travel time uncertainty
has been considered in timetabling in numerous papers, see \cite{LLB18} for a survey.
Transit lines have an impact on the timetable and hence on its robustness.
In order to \emph{evaluate} a line concept with respect to the traveling times of the passengers
in schedule-based transit, a simulation approach can be used that considers not only
the transit lines but also the timetable and the vehicle schedules as follows:
\smallskip

\begin{enumerate}
  \item Based on the line concept, design a timetable and a vehicle schedule (as in
    Figure~\ref{fig-stages}).
  \item Generate a set $\cU$ of (typical) driving times, including unforeseen delays.
  \item Route the passengers in the trip-based aggregation level for every scenario in $\cU$.
  \item Aggregate their traveling times to, e.g., the expected or the worst case
    traveling time.
\end{enumerate}
\smallskip

{\sf LinTim} (see \cite{lintimhp})
is a tool that can be used for such an evaluation as done in \cite{LinTim12}.
A systematic evaluation of a large set of possible disruptions has been suggested in
\cite{FMRSS17,FMRSS18}. For Step 2, the authors generate delays on each single train run,
a travel time increase on each link, and a temporary blocking of each station.
For each of these scenarios the additional delay encountered by the passengers is
computed and summarized in a robustness value. The evaluation has been used within a
machine learning approach for improving the robustness of timetables in 
\cite{MLpaper,mullerhannemann_et_al:OASIcs.ATMOS.2021.3}.

In the following we turn to \emph{analytical} models for dealing with travel time uncertainty.
\smallskip

\textbf{Two-stage robust optimization} 
\cite{CadMar12,CadMar16} use a model for the 
TNDP where passengers are assigned to routes by the planner
in the link-based level of detail and optimize a weighted sum of coverage by transit,
  costs and passengers' traveling times. As uncertainty, a finite set of disruptions
caused, e.g., by
maintenance work is considered in a two-stage program.
Transit lines are constructed 
in the first stage while the passenger routes are adjusted in the second stage
(in the paper this is called  \emph{recovery action}).
The aim is to find the best set of transit lines in the worst case.
\smallskip

\textbf{Two-stage stochastic optimization}
\cite{sadrani2022optimization} assume a log-normal distribution of
driving times between stops. They sample their values from the scenario set within the
optimization. The goal is to minimize the expected travel time of the passengers.
Their routes are determined by the planner on trip-level together with a scenario-based
schedule and under consideration of crowding in the second stage.
The mean and the standard deviation needed are assumed to be given, but analyzed
within a sensitivity analysis.

\cite{Cadarso17} modify the
deterministic model of \cite{CadMar12,CadMar16} by making the driving
times stochastic. They use a finite scenario set $\cU$. 
The aim is to minimize the expected value of coverage by transit as performance
indicator, where the lines are determined in the first stage and the
routes of the passenger are 
adjusted by the planner in the second stage
within the link-based aggregation level.
Due to the finiteness of the uncertainty set the problem can be
formulated as a linear integer program. The results are compared with the
Expected Value model and two different Value-at-Risk formulations.
\medskip

There are also other approaches to increase
robustness of line plans under travel time uncertainty.
Some authors argue that an
equal distribution of passengers through the network may decrease the probability
for delays since crowded arcs tend to be more sensible to delays and since a disruption
will affect less passengers in the worst case if there are no overcrowded arcs.
This has been done in \cite{CadMar12,LapMarMesaOrt11} where constraints on the
passenger routes make sure that demand is forced to disperse on the whole network.
Instead of distributing the passenger routes, \cite{SchSchw13} aim at dispersing the
lines more equally on the network within a game-theoretic model: here,
each line is interpreted as a player, and the players compete for the limited
infrastructure when choosing their frequencies. Both approaches work on the
link-based aggregation level.

\cite{YaoHuLuGaoZhang13} add buffer times to the travel times in the transit line
planning model to account for potential delays 
In \cite{YanLiuMengJiang2013}, a transit travel time reliability constraint is introduced
within a stochastic model requiring that
the probability that the travel time for an OD-pair is larger than a
certain threshold is small.


\section{Extensions and related problems}
\label{sec-extensions}

\subsection{The skip-stop problem}\label{sec-skipstop}
So far, lines have been defined as paths in the public transportation
network PTN, usually by their sequence of links. This determines also
the stops of a line. In the \emph{skip-stop problem} the lines 
are given and fixed. The decision to be made is to determine their stops.
More precisely, for every line the skip-stop problem decides, 
which of the stops along its path should be visited by the line and
which of them can be skipped.

Formally, assume that a line $l$ passes links $e_1=(u,v)$
and $e_2=(v,w)$, i.e., the path for $l$ includes the sequence
\[ l=(\ldots, (u,v), (v,w), \ldots ). \]
If the stop $v$ is \emph{skipped}, the two consecutive links
$(u,v)$ and $(v,w)$ are merged into one link $(u,w)$ and we receive
\[ l'=(\ldots, (u,w), \ldots ) \]
instead. Note that the PTN needs to be transitive when merging edges
which is a realistic assumption since the bus/train can just pass a stop/station
without stopping there.
\medskip

For the model, assume a line $l$ which is
operated every $p$ minutes. The idea is to split the line into
two lines which are operated alternately,
let's say into a red line and into a green line. The red line serves
only red stations $V^{\textup{red}} \subseteq V$ and the green line serves
only green stations $V^{\textup{green}} \subseteq V$. The set of green
stations and the set of red stations does not need to be disjoint.
It is required that every station is served by at least one of the two
lines, i.e., that
\[ V^{\textup{red}} \cup V^{\textup{green}}= V. \]
For every station that is not visited by both lines, driving time is saved
because the line does not need to stop. The saved time can be used to increase
frequency such that more passengers can be transported.
This is a reason to skip as many stops as possible since capacity is
a crucial issue in crowded metropolitan areas. Skip stop problems in this context
have been developed, e.g., in Asia
\cite{CYL14,Cacchiani2020,RHH19,park2013column} or in South America \cite{Munoz13}.
The saved time can also be used as buffer to increase robustness \cite{Cacchiani2020}.
\medskip

For passengers traveling between red stations or between green stations,
the travel time decreases, too. However, passengers who wish to travel between
a red and a green station need an additional transfer. Hence,
skipping too many stations makes it 
inconvenient for the passengers that now do not have direct connections any more
and might even need to travel a short piece into the wrong direction.
\medskip

In the literature, the term skip-stop is used and analyzed from many different perspectives and there exist several approaches of solving the skip-stop problem
considering various objective functions that have to be minimized.
\cite{Vuchic76} analyzes the differences between the skip-stop operation
and the standard operation with respect to
the passengers time savings. Within a parametric approach
\cite{Munoz13} considers a continuum approximation approach
focusing on the density of the stations visited by both lines.
Evenly distributed demands and identical stopping-times are assumed for analyzing
the influence of different parameters of the model, especially the length of the lines.
A multicriteria-optimization model considering waiting times, travel times
and driving times of the trains is developed by \cite{CYL14}.
In a suburban area, \cite{XLHW14} aim at a solution 
with one local line visiting all stations and one express line skipping some
of the stations. In addition to minimizing the total passenger
traveling time and the riding times of the vehicles, the total number of used
trains (for a given number of passengers) is taken into account.
\medskip

Combining the choice of a stopping pattern with the computation of the concrete
timetable has been investigated in \cite{YueWangZhouTongSaat16,JiangCaccToth17},
and robust models hedging against demand uncertainty for the integrated problem are developed in \cite{YinTangYangGaoRan16,Cacchiani2020,JiangCaccToth17}.
These papers use rail operations as examples. Demand uncertainty has also been
investigated for skip-stop bus operations. Here, \cite{MZL20} minimizes the costs taking
into account the uncertain
boarding times of the passengers, while \cite{LYQZ13,LSC14,RHH19} develop robust
strategies for skipping stops if a bus is behind schedule in the deadheading problem.
A quite different application for skip-stop models is presented by \cite{LAKB22} who
propose a skip-stop model that includes social distancing constraints and
enables to plan better schedules while capacities of the vehicles are drastically
reduced (due to the COVID-19 pandemic).

\subsection{Transit network design with seasonal demand}
Which line concept is to be considered optimal highly depends on the demand.
Transit line planning for \emph{uncertain} demand has already been addressed in
Section~\ref{sec-demand-uncertainty}. However,
  demand changes are typical in public transport: at night, there is less demand than
  during the day, in peak-hours we have high-demand periods. We call these different
  demand periods \emph{seasonal demand}. Let us assume that we have
  $K$ different 
demand seasons, each of them characterized by an OD-matrix
  $\OD_j$, $j=1,\ldots,K$.
    When designing transit lines, the two extremes for dealing with seasonal demand
  are the following:
  \begin{itemize}
  \item One might find an optimal line concept $(\cL^j, (f_l)_{l \in \cL^j})$
    for each of the OD-matrices $j=1,\ldots,K$. This would serve each of the demand seasons
    best, but the resulting line concepts might be very different from each other.
  \item The other extreme is to find only \emph{one} line concept $(\cL, (f_l)_{l \in \cL})$
    which should suit all different demand seasons. This is easy to memorize but most likely
    too expensive or too crowded.
  \end{itemize}
  Let us denote $f^j:=(f_l)_{l \in \cL^j}$ the vector containing the frequencies of line concept $\cL^j$.
  
  Seasonal demand and the above two options have already been mentioned in
  in \cite{amiripour2014designing}. 
  The authors suggest to find one line concept $(\cL^r,f^r)$ minimizing the expected
  travel time over all seasons. They call it \emph{robust} since they additionally
  require that it should on average over all seasons be better than any of the
  optimal line concepts $(\cL^j, f^j)$. With $F_j(\cL, f)$ being the travel time of $(\cL,f)$ 
  evaluated for $\OD_j$ they require for all seasons $i$ that
  \[ \sum_{j=1}^k \alpha_j F_j(\cL^r, f^r) \leq \sum_{j=1}^k \alpha_j F_j(\cL^i,f^i)\]
  where $\alpha_j$ is the length of the season with OD-matrix $\OD_j$.
  
  Based on the observations that demand of the metrobüs in Istanbul is extremely
  unsteady (see \cite{borndorfer2018line}) 
  \cite{csahin2020multi} are the first to define a multi-period transit line planning
  problem: they search for a line concept for each of the periods and couple them
  by resource constraints which make sure that the resources flow from one period to the
  next. The model is stated in general form starting from the cost model \eqref{costmodel}
  as basis.
  \medskip
  
  In \cite{ZhuXuWangVan22,DuNiVan22}, the authors take a passenger-oriented point of
  view and evaluate line concepts under
  various OD matrices within a case study. They suggest to find a line concept for the
  morning peak hour (since this is the busiest demand season) and for other
  seasons to update the frequencies while keeping the same line plan.
  The model in \cite{multiperiod} goes a step further. Here, a line concept
  $(\cL^j, f^j)$ is determined for each of the demand seasons, where coupling constraints
  ensure sufficient similarity of the different line concepts:
  \[ \textup{diff}\left( (\cL^j, f^j),(\cL^k, f^k)\right) \leq \alpha \ \ \forall k,l=1,\ldots,K\]

  Different (dis)similarity measures are discussed and tested on the
  basic cost model \eqref{costmodel}: Similarity of line concepts can be evaluated
  by looking at the differences of the frequency vectors (defined on $\Pool$), i.e.,
    $\textup{diff}\left (\cL^j, f^j),(\cL^k, f^k)\right):=\|f^j-f^k\|$
  which is interesting also for the case in which only the frequencies change and the lines
  stay the same for all seasons.
  A simple measure is to investigate only the number of lines that are added
  or that change. The most promising approach starts by defining a 
  dissimilarity $d(l_1,l_2)$ between each pair of lines $l_1,l_2$ first
  and then use these numbers within a transportation problem (similar to the
  Wasserstein-distance) to compute the similarity of the line \emph{concept},
  see \cite{MauUrq09,ZhuXuWangVan22,multiperiod}.

\subsection{Integrating transit line planning with other planning stages}
\label{sec-integration}

The main stages for public transport planning are the design of the
infrastructure, followed by transit line planning, timetabling, vehicle- and
crew scheduling as depicted in Figure~\ref{fig-stages}. Most papers concentrate
on one of these stages.
Proceeding through the stages sequentially, one after another, we receive a
solution which consists of lines, frequencies, a timetable, and vehicle- and crew schedules.
It can be evaluated using similar performance indicators as described in
Section~\ref{sec-performance-indicators}, but the decisions made in later planning stages allow to make more precise estimates, e.g., travel times can be evaluated trip-based instead
of line based and costs for vehicle and crew
can be evaluated based on actual schedules. 
However, proceeding sequentially stage after stage is just a greedy approach and does
in general not find the best possible solution. This has been recognized by many authors.
Solving the planning stages integratedly
would yield an optimal solution but is computationally out of scope. An integer linear
program integrating transit line planning, timetabling, and vehicle scheduling is presented
in \cite{philinediss,SchiSchoe21}
but it can only be solved for tiny networks.
In the following we sketch some work which considers line planning
together with other planning stages.
\medskip

\textbf{Line planning and timetabling.}
Most work on integrating transit line planning with another planning stage concerns
the integration of line planning and timetabling. There are a few exact models (given
as mixed-integer programs), e.g., in \cite{Puerto20} who
model the line selection, frequency setting and timetabling problem respecting
capacities and passenger flows, or in \cite{IrnStei18}
who choose lines and schedules simultaneously for a line corridor in the
context of intercity bus lines. Since such models are hard to solve, different types
of heuristic approaches have been proposed.
\begin{itemize}
\item Iterative approaches switch between transit line planning and timetabling iteratively
  trying to improve the line concept and/or the timetable in every step. Depending on the submodels
  used for transit line planning or timetabling, different improvement strategies exist.
  In \cite{Sch16} lines can be broken at stations and re-arranged within a matching problem
  in order to reduce the number of transfers (but still respecting  the timetable).
 In \cite{BurgBullVanLusby17} the lengths of the lines are changed in each iteration allowing a more
  flexible and hence more robust timetable. Within a railway context
  \cite{YanGov2019combined} adapt the maximal frequencies in each iteration to deal with
  infrastructural limitations while \cite{zhang2022integrated} (also in a railway context)
  adds a coupling constraint whose level of satisfaction is used within a feedback loop in the
  iterative process.
\item A look-ahead strategy is followed by
  \cite{Fuchs22} who aims at finding a line plan that admits a feasible timetable. Within an iterative approach a small set of incompatible services is identified and modified until a feasible
  timetable exists.
\item An evolutionary approach for integrating transit line planning and timetabling is suggested
in \cite{KR13}. 
\end{itemize}
\smallskip

\textbf{Line planning, timetabling, and vehicle scheduling.}
Exact formulations for an integration of these three planning stages are given in
\cite{philinediss,SchiSchoe21} as linear integer programming models using trip-based aggregation
and routing passengers by the planner.
Since it is computationally hard to solve these models, several other models and approaches have
been developed.

\begin{itemize}
\item \cite{liebchenmoehring2007} deal mainly with timetabling integrating vehicle scheduling and
allowing slight rearrangements of lines within their approach.
\item \cite{Sch16} proposes to re-optimize all three stages step by step within an iterative
  approach. As framework, the so-called eigenmodel \cite{Sch16} is developed. Its nodes
  represent subproblems and the edges indicate which of these problems can be solved after another.
  In \cite{JaeSch18} it is shown that the
  approach always converges to a local optimum in a finite number of improvement steps
  for public transit planning. \cite{PSSS17} use this idea by changing
  the order of the subproblems in a sequential approach: they start with vehicle scheduling
  followed by timetabling and line planning
  instead of processing in the classic way indicated in Figure~\ref{fig-stages}.
\item A look-ahead approach is followed by \cite{LiTang23} who develop a line plan
  together with a timetable integratedly and in such a way that it is a feasible input
  for the vehicle scheduling problem. This is modeled
  as a large nonlinear integer optimization problem which is solved by decomposition approaches.
\end{itemize}
\smallskip

\textbf{Integrating transit line planning with other stages.}
The integration of transit line planning and vehicle scheduling (without looking
at a timetable, but just forming the vehicles' routes) has been treated in 
\cite{MicSch07,pternea2015sustainable,SchiSti22}. Also the integration of transit line
planning with tariff planning has been considered, see
\cite{zhou2023joint} or \cite{tirachini2014multimodal} who combine frequency setting and
the determination of the bus sizes with the optimization of fares and congestion tolls.

\subsection{Parametric transit line planning} \label{section-parametric}

There is a stream of literature 
that regards line planning and frequency setting from a \emph{continuum approximation} perspective. Rather than designing a line concept and/or frequencies for a specific city, represented by a PTN and OD-pairs on the PTN, the aim is to make general statements about optimal design of line concepts, frequencies, and capacities.
  An early contribution in this area is \cite{mohring1972optimization}, who, coming from a microeconomics perspective, analytically studied the problem of determining optimal frequencies and the number of bus stops for a so-called  \emph{steady-state route}, where passengers' origins and destinations are uniformly distributed along the route and bus stops are to be spaced evenly along the route. An optimal frequency is obtained by balancing operator cost (which increase with frequency) and passenger waiting times (which decrease with frequency). This leads to the famous \emph{square root formula} which states that   the  frequency should be proportional to the square root of the number of passengers. This line of research has been developed 
further in, e.g., \cite{jansson1980simple,jara2003towards,jara2009effect}. \cite{jara2003single}
are able to generalize the square root result to lines on networks.

Based on these results there is also a stream of publications that consider transit network design on \emph{parametric networks}, i.e., networks whose topology and demand structure are fully determined by a few input parameters. The goal of most of these publications is to determine optimal design parameters for stylized line concepts (\emph{feeder trunk}, \emph{hub and spoke} etc.) and to specify \emph{optimal} line concepts based on the input parameters with respect to a weighted sum of costs and travel time.

\cite{fielbaum2016optimal} represent the area for which the line plan is to be designed as a network consisting of a central business district node (CBD) 
connected to several subcenter nodes which are in turn connected to neighboring subcenters and to periphery nodes. Distances between CBD and subcenters and between subcenters and periphery nodes, as well as the demand structure between the different types of nodes are controlled parameters. 
They analyze the performance of four general network structures (\emph{direct lines}, \emph{feeder trunk}, \emph{hub and spoke}, and \emph{exclusive}) 
 under some common transit network design assumptions.
 For each of the proposed structures, they determine optimal frequencies and capacities for the lines
 with respect to an objective function composed of operational costs (depending on the total number of vehicles needed and the vehicle capacity, but notably not on the distance driven)
 and total passenger travel time (composed of in-vehicle time, adaption time, and transfer penalties).
 Due to the regular network structure, they are able to express the different cost components as functions of the parameters and to compute optimal frequencies and capacities numerically when parameter values are given.
 \cite{fielbaum2021lines} enrich the model from  \cite{fielbaum2016optimal} by \emph{line spacing}. 
\cite{Parametric18} compare the performance of the general line network structures proposed in \cite{fielbaum2016optimal} with line plans found by approaches developed for transit line planning problems with specific (non-parametrically given) PTN and demand. 
\cite{masing2022price} propose a mixed-integer linear program to compute (symmetric) line concepts on parametric cities and discuss the price of symmetry.

Analyses similar to  \cite{fielbaum2016optimal} have been done also for other underlying network structures like grid networks  \cite{daganzo2010structure,chen2015optimal} and radial networks \cite{badia2014competitive,chen2015optimal}.

\subsection{Further related problems}\label{sec-relatedproblems}

During disruptions or construction work at rail infrastructure, the corresponding connections are often bridged by buses, leading to ad-hoc route generation, route selection and frequency assignment problems. The \emph{bus bridging} problem has been studied in, e.g., \cite{kepaptsoglou2010bus, vanderHurk2016shuttle, gu2018plan, deng2018design}, for further references see also the corresponding section in the recent survey on urban rail disruption management by \cite{wang2024urban}.

\emph{On-demand transportation} does often not operate on fixed routes, but
is planned 
  flexibly based on passengers origins and destinations. This transport mode has been used for a long time to transport people with reduced mobility access where requests are normally known in advance.
Recently, it has sparked new interest, also in its online version, with the advent of autonomous vehicles on the horizon. See, e.g., \cite{cordeau2006branch,ropke2009branch,gaul2021solving}. A recent survey is given in \cite{vansteenwegen2022survey}.

  \cite{HoHuSchiSchoUr-atmos} study the \emph{edge (link) investment problem}: which segments of a bus line should be upgraded to be used by bus rapid transport, assuming that there are multiple parties responsible for the investments. 
   
  To avoid \emph{bus bunching}, different control strategies are possible, such as waiting, stop-skipping, short-turning. See, e.g.,  \cite{daganzo2009headway,tian2022short} and references therein. \cite{gkiotsalitis2023exact} investigate the problem for rail-based transport. A survey is given in \cite{rezazada2022public}.

There is also a stream of research dealing with the question how to \emph{visualize line plans} according to criteria like preservation of network topology and relative position of stations, uniform link lengths between stations, few bends along individual metro lines, or limited number of link orientations. 
A mixed integer linear programming approach to the problem has been described in \cite{nollenburg2010drawing}. This approach has been used to generate the line plan in Figure~\ref{fig-lineplan2}. For more on line plan visualization, see \cite{w-dsms-07,DBLP:journals/cgf/WuNTRN20}.

\section{Outlook and further research}
\label{sec-conclusion}

We hope that this article made a small contribution to the following
two goals which we consider as important for a further development
of the field of transit planning. The first is to bring researchers from
different communities even closer together. We believe that cooperations will help
to improve research on optimization of public transport significantly.
The second challenge is to use more results of transit line planning 
for optimizing real-life applications, in particular in view of
climate changes, and eventually make transit line
planning procedures available in practice whenever they are needed.

\section*{Acknowledgments}
The authors thank  Martin Nöllenburg und Alexander Wolff for providing Figure\ref{fig-lineplan2}.

\small
\bibliographystyle{apalike}
\bibliography{eigen,linpla,references,integration}


\end{document}